\documentclass[reqno,12pt]{amsart}
\usepackage[utf8]{inputenc}
\usepackage[left=1.2in, right=1.2in, top=1in]{geometry}

\usepackage{esint}
\usepackage{amsmath,amssymb,amsthm,mathrsfs,color,times,textcomp,verbatim,mathtools}
\allowdisplaybreaks
\usepackage{xcolor}
\usepackage[colorlinks=true]{hyperref}
\hypersetup{urlcolor=blue, citecolor=red, linkcolor=blue}
\usepackage[square,numbers]{natbib}

\numberwithin{equation}{section}
\theoremstyle{plain} 
\newtheorem{theorem}{Theorem}[section]
\newtheorem{proposition}[theorem]{Proposition}
\newtheorem{lemma}[theorem]{Lemma}
\newtheorem*{conjecture}{Conjecture}

\newtheorem{corollary}[theorem]{Corollary}
\newtheorem{claim}{Claim}

\usepackage{graphicx}

\newtheorem{remark}[theorem]{Remark}

\title[Quantitative stability of  harmonic maps]{quantitative stability of  harmonic maps from $\mathbb{R}^2$ to $\mathbb{S}^2$ with a higher degree}
\author{Bin Deng}
\address{School of Mathematics and Statistics, Wuhan University,
Wuhan, Hubei Province, P.R. China, 430072}
\email{dbmath@whu.edu.cn}
\author{Liming Sun}
\address{Academy of Mathematics and Systems Science, the Chinese Academy of Sciences, Beijing 100190, China.}
\email{lmsun@amss.ac.cn}
\author{Jun-cheng Wei}
\address{Department of Mathematics, Chinese University of Hong Kong, Shatin, Hong Kong.}
\email{jcwei@math.ubc.ca}
\date{\today \,(Last Typeset)}

\def\R2{\mathbb{R}^2}
\def\S2{\mathbb{S}^2}

\def\Re{\textup{Re\,}}
\def\im{\textup{i\,}}
\def\CC{\mathbb{C}}

\def\J{\mathcal{J}^r}
\def\I{\mathcal{I}^r}
\def\L{\mathcal{L}}
\def\KK{\mathcal{K}^r}
\def\W{\mathcal{W}_{\Phi}}
\def\H{\mathcal{H}_{\Phi}}
\def\dH{\dot{\mathcal{H}}_{\Phi}}

\begin{document}

\begin{abstract}
For degree $\pm 1$ harmonic maps from $\mathbb{R}^2$ (or $\mathbb{S}^2$) to $\mathbb{S}^2$, Bernand-Mantel, Muratov and Simon \cite{bernand2021quantitative} recently establish a uniformly quantitative stability estimate. Namely, for any map $u:\mathbb{R}^2\to \mathbb{S}^2$ with degree $\pm 1$, the discrepancy of its Dirichlet energy and $4\pi$ can linearly control the $\dot H^1$-difference of $u$ from the set of degree $\pm 1$ harmonic maps. Whether a similar estimate holds for harmonic maps with a higher degree is unknown. In this paper, we prove that a similar quantitative stability result for a higher degree is true only in a local sense. Namely, given a harmonic map, a similar estimate holds if $u$ is already sufficiently near to it (modulo M\"{o}bius transforms) and the bound in general depends on the given harmonic map. More importantly, we thoroughly investigate an example of the degree 2 case, which shows that it fails to have a uniformly quantitative estimate like the degree $\pm 1$ case. This phenomenon shows the striking difference between degree $\pm1$ ones and higher degree ones. Finally, we also conjecture a new uniformly quantitative stability estimate based on our computation.
\end{abstract}

\maketitle

\setlength{\parskip}{0.5em}
\section{Introduction}\label{sec:intro}
\subsection{Motivation and main results} 
The analysis of critical points of conformally invariant Lagrangians has drawn much attention since 1950, due to their important applications in physics and geometry. One of the prominent examples is harmonic maps $u:\R2\to \mathbb{S}^2$, which are critical points of the following Dirichlet energy
\begin{align} 
    \mathcal{E}(u)=\frac{1}{2}\int_{\R2}|\nabla u|^2dx.
\end{align}
It is well-known that the critical points $u$ will satisfy
\begin{align}\label{eq:harmonic}
    \Delta u+|\nabla u|^2u=0,\quad \text{in }\R2.
\end{align}
In the special case of mapping from $\R2$ to $\S2$, all the harmonic maps have been classified. For instance, see 
\cite[11.6]{eells1978report} and \cite[Section 2.2]{helein2008harmonic}. To state the result, we introduce 
\begin{align}\label{constr-p/q-2}
    \mathcal{A}_d=\{(p,q):p/q\text{ is an irreducible rational function of }z,  \max\{\deg p,\deg q\}=d\}
\end{align}
for any $d\in \mathbb{Z}_{\geq 1}$. Here $z$ is the complex variable in $\CC=\R2$. When $d\in\mathbb{Z}_{\leq -1} $, we also introduce the notation  ${\mathcal{A}}_d=\{(\bar p,\bar q):(p,q)\in \mathcal{A}_{|d|}\}$. Throughout this paper, we assume that $d$ is an integer with $|d|\geq 1$.
\begin{proposition}\label{prop:class}
A map $u$ from $\R2(=\CC)$ to $\S2$ is harmonic if and only if $u$ is holomorphic or anti-holomorphic. More precisely, if $\deg(u)=d $ then $u =\mathcal{S}(p/q)$ 
where $(p,q)\in \mathcal{A}_d$. Here $\mathcal{S}$ is the stereographic projection from $\CC\to \S2\setminus\{N\}$ (see \eqref{def:stereo-proj}).
\end{proposition}

For any pair $(p,q)\in \mathcal{A}_d$, we can normalize $p$ to be a monic polynomial such that $p/q$ stays the same. We call such pair $(p,q)$ to be canonical. Thus we define
\begin{align}
    \mathcal{A}_d^m=\{(p,q)\in \mathcal{A}_d:p\text{ is monic}\}.
\end{align}
Then for any harmonic map $\Phi:\R2\to \S2$ it can be represented by a unique canonical pair $(p,q)\in \mathcal{A}_d^m$.

It is also known that harmonic maps from $\R2$ to $\S2$ achieve minimal Dirichlet energy within its homotopy class by the work of \citet{lemaire1978applications} and \citet{wood1974} (also see \cite[(11.5)]{eells1978report}). 
\begin{lemma}\label{intro:lem:minimizer}
Suppose that $u:\R2 \to \S2$ with $\mathcal{E}(u)<\infty$. Then $\mathcal{E}(u)\geq 4\pi|\deg (u)|$. The equality holds if and only if $u$ is harmonic.
\end{lemma}

A natural question is whether the discrepancy $\mathcal{E}(u)-4\pi |\deg(u)|$ can quantitatively control the difference of $u$ from the harmonic maps. Such a type of question has been raised for many other topics. For instance, \citet{brezis1985sobolev} ask a similar question to the classical Sobolev inequality on $\mathbb{R}^n$. Later \citet{bianchi1991note} obtain a quantitative stability estimate in the spirit of  \eqref{stab-1}. (See also recent work of \citet{FigalliZheng}.)   \citet{fusco2008sharp} prove sharp quantitative stability about isoperimetric inequality.

Recently \citet*{bernand2021quantitative} prove quantitative stability for degree $\pm 1$ harmonic maps from $\R2$ to $\S2$ as reformulated in the following. The proof in \cite{bernand2021quantitative} has been simplified by \citet{hirsch2021note} and \citet{topping2020rigidity}.

\begin{proposition}
There exists a universal constant $C$ such that for every $u:\R2\to \S2$ with $\mathcal{E}(u)<\infty$ and $\deg(u)=1$, there exist $a,b,c,d\in \mathbb{C}$ with $ad-bc\neq 0$ such that 
\begin{align}\label{stab-1}
    \int_{\R2}\left|\nabla u-\nabla \mathcal{S}\left(\frac{az+b}{cz+d}\right)\right|^2\leq C(\mathcal{E}(u)-4\pi).
\end{align}
If $\deg(u)=-1$, then the above statement holds for $\bar z$.
\end{proposition}

The above theorem leaves us with an intriguing question for harmonic maps with a higher degree. We have addressed a similar question for half-harmonic maps in \cite{deng2021non}. There we have shown that similar quantitative stability for higher degree ones is only true in the local sense. In this paper, we shall prove that a similar phenomenon happens here. More precisely, given a harmonic map (or a compact set of harmonic maps), there is a local stability result near it. The bound in general will depend on the given harmonic map (or the compact set). 

To that end, let us introduce some notations. 
For any two complex polynomials $p, \tilde p$ on $z$ (or $\bar z$), we define $|p-\tilde p|_\infty$ to be the maximum of all the modulus of each coefficient of $p-\tilde p$. Thus $\mathcal{A}_d^m$ becomes a metric space equipped with the distance $|\cdot |_\infty$.
When we call $\Omega$ a compact set of degree $d$ harmonic maps from $\R2$ to $\S2$, it actually means  there exists a compact set of $\mathcal{A}_d^m$, say $\mathcal{A}_{\Omega}^m$, such that $\Omega=\mathcal{S}(\mathcal{A}_{\Omega}^m)$. Here compact is in the sense of $|\cdot|_\infty$-topology. Namely, for any sequence of $(p_k,q_k)\in \mathcal{A}_\Omega^m$, it has a subsequence $(p_{k'},q_{k'})$ and some $(p_*,q_*)\in \mathcal{A}_\Omega^m$ such that $|p_{k'}-p_*|_{\infty}+|q_{k'}-q_*|_{\infty}\to 0$ as $k'\to \infty$. For any map $u:\R2\to \S2$, we define the $\dot H^1$-distance in the following way:
\begin{align} 
    dist_{\dot H^1}(u,\Omega)=\inf_{\Phi\in \Omega}\left(\int_{\R2}|\nabla u-\nabla \Phi|^2\right)^{\frac12}.
\end{align}

\begin{theorem}[\textbf{local stability}]\label{thm:main-2}
Suppose $\Omega$ is a compact set of degree $d$ harmonic maps. There exist two constants $\eta=\eta({\Omega})$ and $C=C(\Omega)$ such that if $u:\R2\to \S2$ with $dist_{\dot H^1}(u, \Omega\circ F)<\eta$ for some M\"obius transformation $F:\R2\to \R2$, then there exists a harmonic map $\Phi$ which makes the following hold
\begin{align}
    \int_{\R2}|\nabla u-\nabla \Phi|^2\leq C(\mathcal{E}(u)-4\pi|\deg(u)|).
\end{align}
Moreover, $\Phi\circ F^{-1}$ is near to $\Omega$ in the sense of $|\cdot|_\infty$-topology.
\end{theorem}

For any harmonic map, Proposition \ref{prop:class} says that there is a family of harmonic ones near it. We take $\Omega$ to be a compact neighborhood of a given harmonic map, then the above theorem indicates that there is a local version of stability.
Similar to degree $\pm 1$ harmonic maps, one attempts to remove the dependence on $\Omega$ and get a \textbf{uniformly} quantitative stability estimate. However, the following theorem says that this is not true.

\begin{theorem}[\textbf{non-uniform stability}]\label{thm:main-1}
For any large constant $M>0$, we can find $u:\R2\to \S2$ with $\deg(u)=2$ such that for any $(p,q)\in \mathcal{A}_2$  
\begin{align}\label{stab-high}
    \int_{\R2}|\nabla(u-\mathcal{S}(p/q))|^2> M \left(\mathcal{E}(u)-8\pi\right).
\end{align}
\end{theorem}

The above theorem indicates that there are some fundamental differences between the case of degree $\pm1$ and higher degree. To name one, for degree $1$ (resp. $-1$) harmonic maps, one can compose it with a certain M\"obius transformation of $\R2$ such that it becomes the stereographic projection $\mathcal{S}:\R2\to \S2$ (resp. $\bar{\mathcal{S}}$, i.e. the conjugate of $\mathcal{S}$). Since both sides of \eqref{stab-1} is invariant under M\"obius transformations, essentially \eqref{stab-1} is a quantitative stability near $\mathcal{S}$ (resp.  $\bar{\mathcal{S}}$).  However, for higher degree ones, one can not transform a harmonic map to an arbitrary another one by M\"obius transformations. 


Of course, one wishes to have uniformly quantitative stability like \eqref{stab-1} for higher degree harmonic maps. Since Theorem \ref{thm:main-1} says that this is not possible with the naive extension, then one probably needs to minus more things in the square on the left-hand side of \eqref{stab-high}, or strengthen the right-hand side to some nonlinear expression of $\mathcal{E}(u)-4\pi |\deg (u)|$. Actually,  we make the following conjecture in the higher degree case:
\begin{conjecture} Let $|d|\geq 2$.  There exists a universal constant $C=C(d)$ such that for every $u:\R2\to \S2$ with $\mathcal{E}(u)<\infty$ and $\deg(u)=d$, there exists $(p,q)\in \mathcal{A}_{d}$ such that 
\begin{align}\label{stab-high-uni}
        \int_{\R2}|\nabla(u-\mathcal{S}(p/q))|^2\leq C \left(\mathcal{E}(u)-4\pi|d|\right)\left(1+\big|\log (\mathcal{E}(u)-4\pi|d|)\big|\right).
\end{align}
\end{conjecture}
This conjecture is based on the explicit computation of the example we construct in the proof of Theorem \ref{thm:main-1}. In that example, one gets $r^{-4}$ on the left hand side of \eqref{stab-high} and $r^{-4}|\log r|^{-1}$ on its right hand side for some $r\gg 1$ (see \eqref{RHS-log} and \eqref{LHS-4}). Plugging these two facts into \eqref{stab-high-uni}, they exactly make the two sides comparable. 

\begin{remark}
    After our paper appeared on arXiv, there is an interesting preprint addressing the quantitative stability for maps from $\mathbb{S}^2$ to $\mathbb{S}^2$ with general degree by \cite{rupflin2023sharp}.
\end{remark}
To mention a few related works, the sharp quantitative stability of the Euler-Lagrange equation of Sobolev inequality on $\mathbb{R}^n$ also varies according to the number of bubbles and the dimension $n$, readers can consult \cite{ciraolo2018quantitative,figalli2020sharp,deng2021sharp}. In a different direction from ours, \citet{topping2004repulsion} uses the torsion of an almost-harmonic map $u$ to control its $\mathcal{E}(u)-4\pi |\deg(u)|$, which can be considered as another notion of quantitative stability.

\subsection{Comment on proofs}
The locally quantitative stability theorem is proved by using the non-degeneracy result of the harmonic maps (see \cite{gulliver1989rate,chen2018nondegeneracy}). The non-degeneracy (also called integrability) of the linearized operator implies that it has a spectral gap, which can be used to prove local stability near one harmonic map. This is also how \cite{bernand2021quantitative} proves \eqref{stab-1} for $\deg(u)=\pm 1$. We generalize the approach of \cite{bernand2021quantitative} to a higher degree case.

The proof of Theorem \ref{thm:main-1} follows the general framework as the one in our work of half-harmonic maps \cite{deng2021non} with new essential difficulties. First, the Jacobian of the kernels is uniformly non-degenerate in \cite{deng2021non}, while it degenerates as the parameter goes to infinity in the case of harmonic maps (cf. \eqref{det-Jr}). One needs to expand up to the third order in the computation to observe this fact, which makes the process substantially more involved. Second, the trick of choosing a vector field corresponding to the rotation in \cite{deng2021non} does not work here, which is the heart of the construction. Fortunately, we leverage the degenerate tendency of the Jacobian to find a new one (cf. \eqref{def:KK}). We have not gotten a satisfactory explanation about why such a vector field works. 

\subsection{Structure of the paper} In the section \ref{sec:pre}, we give some preliminary of harmonic maps from $\R2$ to $\S2$. In section \ref{sec:stability}, we prove the local stability result. In the section \ref{sec:counterexample}, we construct an example such that \eqref{stab-high} holds, thus Theorem \ref{thm:main-1} is proved. In section \ref{sec:cal-Jac}, we provide some computations which are needed in the previous section.

\section{Preliminary}\label{sec:pre}
Topological degree of a $C^{1}$ map $u$ from $\mathbb{R}^{2}$ to $\mathbb{S}^{2}$ can be defined by the de Rham approach
\begin{align}\label{W-10}
\operatorname{deg}(u)=\frac{1}{4 \pi} \int_{\mathbb{R}^{2}} u \cdot\left(u_{y} \times u_{x}\right)=\frac{1}{4 \pi} \int_{\mathbb{R}^{2}} u_{x} \cdot\left(u \times u_{y}\right).
\end{align}
It is well-known that \eqref{W-10} is equivalent to the Brouwer's degree for all $C^1$ maps. See, for instance, \cite[Chapter III]{outerelo2009mapping}. It is easy to know that the degree is continuous in $\dot H^1(\R2;\S2)$ topology. That is, given any $u\in \dot H^1(\R2; \S2)$, there exists $\eta_1(u)$ such that if $v\in \dot H^1(\R2; \S2)$ with 
\begin{align}\label{H1-contin}
    \|u-v\|_{\dot H^1}\leq \eta_1(u),
\end{align}
then $\deg(v)=\deg(u)$. Moreover, $\eta_1(u)$ can be made uniform to $u$ if $u\in \Omega$ which is a compact set of degree $d$ harmonic maps. In this case, one can replace $\eta_1(u)$ by $\eta_1(\Omega)$.

As discussed in the introduction, the following lemma is already known. For the reader's convenience, we provide direct proof here.

\begin{lemma}\label{lem:minimizer}
Suppose that $u:\R2 \to \S2$ with $\mathcal{E}(u)<\infty$. Then $\mathcal{E}(u)\geq 4\pi|\deg (u)|$. The equality holds if and only if $u$ is harmonic.
\end{lemma}
\begin{proof}
By completing the square it is easy to establish the following identity
\begin{align}\label{nab_u-deg}
    |\nabla u|^2\pm 2u\cdot( u_x\times u_y)=|u_x\mp u\times  u_y|^2.
\end{align}
Integrating on both sides and noticing \eqref{W-10}, one obtains $\mathcal{E}(u)\geq 4\pi |\deg(u)|$. Since Brouwer's degree is invariant by homotopy deformations, we find that $u$ is a minimizer in its degree class. If $\mathcal{E}(u)=4\pi |\deg(u)|$ then $u$ will be a critical point of $\mathcal{E}(u)$, thus it is harmonic.

Conversely, suppose that $u$ is a harmonic map from $\R2$ to $\mathbb{S}^2$. Recall the Hopf differential 
\begin{align} 
    \mathcal{H}(z)=u_y\cdot u_y-u_x\cdot u_x+2\im u_x\cdot u_y.
\end{align}
It is well-known that $\mathcal{H}$ is a holomorphic or anti-holomorphic function on $\mathbb{C}$ (one can use \eqref{eq:harmonic} to verify this directly). Since $\mathcal{E}(u)<\infty$, then $\mathcal{H}(z)\in L^1(\mathbb{C})$. It is easy to show that $\mathcal{H}(z)\equiv 0$. Thus $u_x\cdot u_y=0$ and $|u_x|=|u_y|$. Since $|u|=1$, then $u\perp u_x$ and $u\perp u_y$. Combining these facts, we must have $u_x=u\times u_y$ or $u_x=-u\times u_y$. In any case, it holds that $\mathcal{E}(u)=4\pi |\deg (u)|$.
\end{proof}

Let $z=(x,y)\in\R2=\mathbb{C}$ and $s=(s_1,s_2,s_3) \in\S2$. Define the  stereographic projection  
\begin{align}\label{def:stereo-proj}
    \mathcal{S}:\mathbb{C}\to\S2\setminus\{N\}\quad\text{by}\ s_1=\frac{2x}{1+|z|^2}, s_2=\frac{2y}{1+|z|^2}\ \text{and}\ s_3=\frac{|z|^2-1}{1+|z|^2}.
\end{align}
Alternatively, in complex variable form,
\begin{align}
   \mathcal{S}=\frac{1}{1+|z|^{2}}(2z,|z|^2-1).
\end{align}
Suppose $u=\mathcal{S}(\psi)$ where $\psi$ is a  meromorphic function on $\mathbb{C}$. Then we have
\begin{align}\label{pu}
    \begin{split}
        \partial u =\frac{1}{(1+|\psi|^2)^2}\left( 2(1+|\psi|^2)\partial \psi- 2\psi\partial|\psi|^2, 2\partial|\psi|^2\right).
    \end{split}
\end{align}
Here $\partial$ could be $\partial_x,\partial_y$ or with respect to a real parameter on which $\psi$ depends. It is easy to see that 
\begin{align}\label{pu-norm}
    |\partial u|^2=\frac{4|\partial \psi|^2}{(1+|\psi|^2)^2}.
\end{align}
In particular, for a harmonic map, i.e. $\psi=p/q$ where $(p,q)\in \mathcal{A}_d$, one has
\begin{align}\label{nablaS2}
    |\nabla \mathcal{S}(p/q)|^2=\frac{4(|\partial_x(p/q)|^2+|\partial_y(p/q)|^2)}{(1+|p/q|^2)^2}=\frac{4|\partial_x p q-p\partial_x q|^2+4|\partial_y pq-p\partial_y q|^2}{(|p|^2+|q|^2)^2}.
\end{align}
Suppose that $(p,q)$ satisfies \eqref{constr-p/q-2} with $|d|\geq 1$. Then there exists a constant $C(p,q)$ such that  
\begin{align}\label{decay-4}
    |\nabla \mathcal{S}(p/q)|^2\leq C(p,q)(1+|z|)^{-4},\quad \forall\,z\in \mathbb{C}.
\end{align}  
Moreover,  if $(p,q)$ belongs to a compact set $\mathcal{A}_{\Omega}$ of $\mathcal{A}_d$, then $C(p,q)$ can be replaced by some uniform constant $C(\Omega)$.

The linear equation of \eqref{eq:harmonic} is
\begin{align}\label{eq:linear}
    \L[u](v):=\Delta v+2(\nabla u:\nabla v)u+|\nabla u|^2v
\end{align}
where 
\[\nabla u:\nabla v=\sum_{i=1}^3\nabla u^i\cdot\nabla v^i.\]
We shall use this notation throughout this paper.


Take any harmonic map $\mathcal{S}(p/q)$ where $(p,q)\in \mathcal{A}_d$. Apparently, changing the coefficients of $p$ or $q$ continuously yields a family of harmonic maps. Therefore it generates kernel maps of the linearized operator $\L[\mathcal{S}(p/q)]$. For instance, using \eqref{pu} and taking the derivative with respect to the real (or imaginary) part of a coefficient of $p$ or $q$ will produce a kernel map of $\L[\mathcal{S}(p/q)]$. Conversely, this is also true. It is called  the non-degeneracy (or integrability) of harmonic maps.
\begin{proposition}[\cite{gulliver1989rate,chen2018nondegeneracy}]\label{prop:non-deg}
Suppose $\Phi:\R2\to \S2$ is a harmonic map of degree $d$. Then all bounded maps in the $\ker \L[\Phi]$ are generated by harmonic maps near $\Phi$. In particular, $\dim_{\mathbb{R}}\ker \L[\Phi]=4|d|+2$. 
\end{proposition}
\begin{remark}\label{rem:K}
 Furthermore, if $\Phi=\mathcal{S}(p/q)$, using \eqref{pu}, one can see that each kernel map $K\in \ker\L[\Phi]$ is smooth and depends on $p,q$ smoothly. In addition, \eqref{pu-norm} implies that
\begin{align}\label{decay-K}
|\nabla^jK|(z)\leq C_{K,j}|z|^{-|d|-j}.
\end{align}
\end{remark}

We have the expansion of Dirichlet energy in the following. In the smooth setting, one can compare it with the second variation of Dirichlet energy in \cite[page 7]{lin2008analysis}.
\begin{lemma}\label{lem:2nd-vari}
    Suppose that $\Phi:\R2\to \S2$ is a harmonic map. Assume  $v\in \dot H^1(\R2;\mathbb{R}^3)\cap L^\infty(\R2;\mathbb{R}^3)$ with $v\cdot\Phi=0$. Then for $\varepsilon>0$ small
    \begin{align}\label{lin-2nd} 
        \mathcal{E}(\sqrt{1-\varepsilon^2 |v|^2}\Phi+\varepsilon v)=\mathcal{E}(\Phi)+\frac12\varepsilon^2\int_{\R2}|\nabla v|^2-|\nabla\Phi|^2|v|^2+O_{v,\Phi}(\varepsilon^3)
    \end{align}
    where $|O_{v,\Phi}(\varepsilon^3)|\leq C(|v|_{L^\infty},|v|_{\dot H^1},\Phi)\varepsilon^3$ as $\varepsilon\to 0$.
\end{lemma}
\begin{proof}
We shall choose $\varepsilon$ small so that $\varepsilon|v|<\frac12$. Denote $f=\sqrt{1-\varepsilon^2|v|^2}$.
Then
\[\partial_{\alpha}(\varepsilon v^j+f\Phi^j)=\varepsilon \partial_\alpha v^j+f\partial_\alpha \Phi^j+\partial_\alpha f\Phi^j\]
for $\alpha\in\{x,y\}$ and $j\in\{1,2,3\}$. Thus
\begin{align}
    |\nabla (\varepsilon v+f\Phi)|^2=\varepsilon^2|\nabla v|^2+f^2|\nabla\Phi|^2+|\nabla  f|^2+2\varepsilon f \nabla v:\nabla \Phi+2\varepsilon \partial_\alpha v^j\partial^\alpha f\Phi_j.
\end{align}
Here we have used the Einstein summation convention and $\partial_\alpha \Phi^j\Phi_j=0$ for $\alpha\in \{x,y\}$. We shall integrate the above equation and estimate them one by one. First, note that
\begin{align}
    \int_{\R2} f^2|\nabla \Phi|^2=\int_{\R2}|\nabla\Phi|^2-\varepsilon^2\int_{\R2}|\nabla \Phi|^2|v|^2.
\end{align}
Second, 
\begin{align}
    \int_{\R2}|\nabla f|^2\leq \varepsilon^4\int_{\R2}\frac{|v|^2|\nabla v|^2}{1-\varepsilon^2|v|^2}\leq C(|v|_{L^\infty},|v|_{\dot H^1})\varepsilon^4.
\end{align}
Third, since \eqref{decay-4}, we can integrate by parts to get
\begin{align}
    \int_{\R2}f\nabla v:\nabla \Phi=-\int_{\R2}v:\nabla (f\nabla \Phi)=-\int_{\R2}fv\cdot\Delta \Phi+v^j\partial_\alpha f\partial^\alpha \Phi_j.
\end{align}
The first term on the right-hand side is zero because $\Delta\Phi=-|\nabla\Phi|^2\Phi$. For the second one, we apply H\"older's inequality to get
\begin{align}
    \left|\int_{\R2}v^j\partial_\alpha f\partial^\alpha \Phi_j\right|\leq 
    \left(\int_{\R2}|\nabla f|^2\right)^{\frac12}\left(\int_{\R2}|\nabla \Phi|^2\right)^{\frac12}|v|_{L^\infty}\leq C(|v|_{L^\infty})\varepsilon^2.
\end{align}
Fourth, since $\partial_\alpha v^j\Phi_j=-\partial_\alpha \Phi_jv^j$, we have
\begin{align}
\begin{split}
    \left| \int_{\R2}\partial_\alpha v^j\partial^\alpha f \Phi_j\right|= \left| \int_{\R2}v^j\partial^\alpha f \partial_\alpha \Phi_j\right|\leq C(|v|_{L^\infty})\varepsilon^2.
\end{split}
\end{align}
Combining the above four points, we obtain the conclusion.
\end{proof}

For any $u$ with $\deg(u)=\pm1$, \citet[Lemma 4.3]{bernand2021quantitative} prove that the $\dot H^1$-distance of $u$ to the set of degree $\pm1$ harmonic maps can be achieved. For maps with a higher degree, we show that this is also true when $u$ is already near to them. The following proof is pointed out to us by the anonymous referee, which is simpler than our original ones.

\begin{lemma}\label{lem:inf-ach}
    For any map $u:\R2\to \S2$ of degree $d$, if 
    \begin{align}\label{dist<8pi}
        dist_{\dot H^1}(u,\mathcal{S}(\mathcal{A}_d))=\inf_{\Phi\in \mathcal{S}(\mathcal{A}_d)}\left(\int_{\R2}|\nabla u-\nabla \Phi|^2\right)^{\frac12}<\sqrt{8\pi},
    \end{align}
    then the infimum can be achieved. Furthermore, the set of all minimizers is a compact set in $\dot H^1(\R2,\S2)$.
\end{lemma}
\begin{proof}
Suppose $\Phi_k\in \mathcal{S}(\mathcal{A}_d)$ be a sequence of harmonic maps that is a minimizing sequence for $dist_{\dot H^1}(u,\mathcal{S}(\mathcal{A}_d))$. Since $\int_{\R2}|\nabla \Phi_k|^2dx=8\pi|d|$ for any $k$, then there exists a harmonic map $\Phi_\infty\in C^\infty\cap \dot H^1(\R2,\S2)$ such that $\Phi_k
\rightharpoonup\Phi_\infty$ weakly in $\dot H^1(\R2)$. Let $d_\infty$ be the degree of $\Phi_\infty$. By the lower semicontinuity, one has 
\begin{align*}
    8\pi|d_\infty|=\int_{\R2}|\nabla \Phi_\infty|^2\leq 8\pi |d|.
\end{align*}
Therefore $|d_\infty|\leq |d|$. Notice that if $|d_\infty|=|d|$, then there is no energy loss between $\Phi_k$ and $\Phi_\infty$ hence $\Phi_k$ converges to $\Phi_\infty$  strongly in $\dot H^1(\R2)$ and $\deg(\Phi_\infty)=d$. In this case, $dist_{\dot H^1}(u,\mathcal{S}(\mathcal{A}_d))$ can be achieved by $\Phi_\infty$.

If $|d_\infty|<|d|-1$, then using $\Phi_\infty$ is a harmonic map,
\begin{align*}
dist_{\dot H^1}(u,\mathcal{S}(\mathcal{A}_d))^2 & \geq \int_{\mathbb{R}^2}\left|\nabla u-\nabla \Phi_{\infty}\right|^2 d x \\
& =\int_{\mathbb{R}^2}|\nabla u|^2 d x+\int_{\mathbb{R}^2}\left|\nabla \Phi_{\infty}\right|^2 d x-2 \int_{\mathbb{R}^2} \nabla u \cdot \nabla \Phi_{\infty} d x \\
& =\int_{\mathbb{R}^2}|\nabla u|^2 d x+\int_{\mathbb{R}^2}\left|\nabla \Phi_{\infty}\right|^2\left(1-2 u \cdot \Phi_{\infty}\right) d x 
\end{align*}
Since $|u|=1$, then completing the square implies that $1-2u\cdot \Phi_\infty=|u-\Phi_\infty|^2-1$. Continuing the above inequalities
\begin{align*}
dist_{\dot H^1}(u,\mathcal{S}(\mathcal{A}_d))^2& \geq \int_{\mathbb{R}^2}|\nabla u|^2 d x-\int_{\mathbb{R}^2}\left|\nabla \Phi_{\infty}\right|^2+\int_{\mathbb{R}^2}\left|\nabla \Phi_{\infty}\right|^2\left|u-\Phi_{\infty}\right|^2 d x \\
& \geq \int_{\mathbb{R}^2}|\nabla u|^2 d x-\int_{\mathbb{R}^2}\left|\nabla \Phi_{\infty}\right|^2 d x \geq 8 \pi\left(|d|-\left|d_{\infty}\right|\right)  \geq 8 \pi .
\end{align*}
This contradicts the assumption \eqref{dist<8pi}. Thus $d_\infty=d$ and $\Phi_\infty$ achieves the minimum of $dist_{\dot H^1}(u,\mathcal{S}(\mathcal{A}_d))$.

To prove the compactness of the minimizer, one can start with a sequence of minimizers $\Phi_n$ and repeat the above proof to show $\Phi_n$ converges strongly to some $\Phi_\infty$, which is a minimizer. 

\end{proof}

Recall that any harmonic map $\Phi:\R2\to \S2$ can be represented by a unique canonical pair $(p,q)\in \mathcal{A}_d^m$. We can prove that such representation is continuous from $\|\cdot\|_{\dot{H}^1}$-topology of $\mathcal{S}(\mathcal{A}_d)$ to $|\cdot|_\infty$-topology of $\mathcal{A}_d^m$.


\begin{lemma}\label{lem:pq-near}
    Fix any harmonic map $\Phi=\mathcal{S}(p/q)$ with $(p,q)\in \mathcal{A}_d^m$ being canonical. Then for any $\varepsilon>0$, there is $\eta_2(\Phi,\varepsilon)$ such that if $\tilde \Phi$ is another harmonic map with canonical $(\tilde p,\tilde q)$ satisfying $\|\Phi-\tilde \Phi\|_{\dot H^1}\leq \eta_2(\Phi,\varepsilon)$. Then 
    \begin{align}
        |\tilde p-p|_{\infty}+|\tilde q-q|_{\infty}<\varepsilon.
    \end{align}
    Moreover, if $\Omega$ is a compact set of harmonic maps, then $\eta_2$ can be uniform with respect to $\Phi\in \Omega$.
\end{lemma}
\begin{proof}
    If $d=0$, then $\mathcal{A}_0^m=\{(1,q):q \text{ is a non-zero constant}\}$. It is easy to see the conclusion holds. In the following, we assume $|d|\geq 1$.

    Argue by contradiction. Suppose there exists $\varepsilon_0$ such that for any $k\geq 1$ there exists $\Phi_k=\mathcal{S}(p_k/q_k)$ with canonical $(p_k,q_k)\in \mathcal{A}_d^m$ satisfying $| p_k-p|_\infty+| q_k-q|_\infty\geq \varepsilon_0$ and $\|\Phi_k-\Phi\|_{\dot H^1}<1/k$. 

    
    Since $p_k,q_k$ are all complex polynomials, if all the coefficients of $p_k,q_k$ are all uniformly bounded, then there must exist a subsequence $(p_{k'},q_{k'})$ and a canonical pair $(p_*,q_*)\in \mathcal{A}_d^m$ such that $|p_{k'}-p_*|_\infty+|q_{k'}-q_*|_\infty\to 0$. Letting $k'\to \infty$ in $\|\Phi_{k'}-\Phi\|_{\dot H^1}<1/k'$, one obtains $\Phi_*=\Phi$ and consequently $p_*=p$ and $q_*=q$. This is a contradiction. 
    
    Therefore the coefficients of $p_k,q_k$ can not be all uniformly bounded. We have three possible consequences, namely there exists a subsequence (which we still denote $p_k,q_k$) such that $|p_k/q_k|\to \infty$ a.e. as $k\to \infty$, or a subsequence $|p_k/q_k|\to 0$ a.e. as $k\to \infty$, or a subsequence and a rational function $p_*/q_*$ such that  $p_k/q_k\to p_*/q_*$ locally uniformly on $\R2\setminus\{\text{zeros of }q_*\}$. Moreover, $\max\{\deg p_*,\deg q_*\}<|d|$.
    
    In the first case, $\mathcal{S}(p_k/q_k)\to (0,0,1)$ for almost every $z$. For any $v\in C^\infty_c(\R2;\mathbb{R}^3)$, we have
\begin{align}
    \int_{\R2}\nabla \mathcal{S}(p_k/q_k):\nabla v=-\int_{\R2}\mathcal{S}(p_k/q_k):\Delta v\to 0,\quad \text{as }k\to \infty.
\end{align}
Note that $\|\nabla \mathcal{S}(p_k/q_k)\|_{L^2(\R2;\mathbb{R}^3)}^2=8\pi |d|$ is uniformly bounded. Thus it has a subsequence that weakly converges to $0$ in $L^2(\R2;\mathbb{R}^3)$. However, this leads to 
\begin{align}\label{nab_u-weak}
    \|\nabla \Phi\|_{L^2}=\liminf_{k\to \infty}\|\nabla(\Phi-\mathcal{S}(p_k/q_k))\|_{L^2}=0
\end{align}
which contradicts to the fact that $\mathcal{E}(\Phi)\geq 4\pi |d|\geq 4\pi$.

    In the second cases, we have $\nabla \Phi_k\to 0$ weakly in $L^2(\R2;\mathbb{R}^3)$. Then using the trick in \eqref{nab_u-weak}, we also obtain a contradiction.
    
    In the third case, we must have $|\deg(\mathcal{S}(p_*/q_*))|<|d|$. It can be excluded as before by making $\eta_2$ small enough.
\end{proof}
The following corollary shows that the minimizers of \eqref{inf-dist} should be near to $\Omega$.
\begin{corollary}\label{cor:pq-near}
Suppose $\Omega$ is a compact set of degree $d$ harmonic maps. Then for any $\varepsilon>0$, there exists a constant $\eta_3(\Omega,\varepsilon)$ such that if $u:\R2\to \S2$ with  $dist_{\dot H^1}(u,\Omega)<\eta_3(\Omega,\varepsilon)$, then the following infimum can be achieved
\begin{align}\label{inf-dist}
     dist_{\dot H^1}(u,\mathcal{S}(\mathcal{A}_d))^2=\inf_{\forall(p,q)\in\mathcal{A}_d}\int_{\R2}|\nabla u-\nabla \mathcal{S}(p/q)|^2=\int_{\R2}|\nabla u-\nabla  \Phi|^2.
\end{align}
Moreover, any minimizer $\tilde \Phi$ of \eqref{inf-dist} has canonical representation $(\tilde p,\tilde q)\in\mathcal{A}_d^m$ satisfying 
\[|\tilde p-p|_\infty+|\tilde q-q|_\infty<\varepsilon\]
for some canonical pair $(p,q)\in\mathcal{A}_d^m$ such that $\mathcal{S}(p/q)\in \Omega$. 
\end{corollary}
\begin{proof}
Given any $\varepsilon>0$, we take $\eta_3(\Omega,\varepsilon)=\min\{\eta_1(\Omega),\frac13\eta_2(\Phi,\varepsilon),\sqrt{8\pi}\}$. Then $\deg u=d$.
It follows from Lemma \ref{lem:inf-ach} that there exists a minimizer of \eqref{inf-dist}.
For any minimizer $\tilde \Phi$ of the infimum, we have
\begin{align}
    \|\tilde \Phi-\Phi\|_{\dot H^1}\leq \|\tilde \Phi-u\|_{\dot H^1}+\|u- \Phi\|_{\dot H^1}\leq  2\eta_3<\eta_2.
\end{align}
Thus we can apply Lemma \ref{lem:pq-near} to get the conclusion.
\end{proof}

\section{Local stability}\label{sec:stability}



In this section, we shall prove the local stability result, i.e.\,Theorem \ref{thm:main-2}. Throughout this section, we will assume $\Omega$ is a compact set of degree $d$ harmonic maps from $\R2$ to $\S2$. Moreover, $\Phi$ always denotes a harmonic map.
 Define 
\begin{align*}
    \W(\R2)&=\{v\in L_{loc}^1(\R2;\mathbb{R}^3):\int_{\R2}|v|^2|\nabla \Phi|^2<\infty,v\cdot \Phi=0\},\\
    \H(\R2)&=\{v\in H^1_{loc}(\R2;\mathbb{R}^3):\int_{\R2}|\nabla v|^2+|v|^2|\nabla \Phi|^2<\infty, v\cdot \Phi=0\},\\
    \dH(\R2)&=\{v\in H_{loc}^1(\R2;\mathbb{R}^3):\int_{\R2}|\nabla v|^2<\infty, v\cdot \Phi=0\}.
\end{align*}
Then $\W$ is a Hilbert space with inner product $(v_1,v_2)_{\W}=\int_{\R2}v_1\cdot v_2|\nabla \Phi|^2$; $\H$ is a Hilbert space with inner product $(v_1,v_2)_{\H}=\int_{\R2}\nabla v_1: \nabla v_2+|\nabla \Phi|^2v_1\cdot v_2$.
Similar to Theorem A.1 and A.2 in  \cite{figalli2020sharp}, we can prove the following two lemmas.
\begin{lemma}\label{lem:cpt}
There exists a constant $C(\Phi)$ such that \[\|v\|_{\W(\R2)}\leq C(\Phi)\|v\|_{\dH(\R2)}.\]
Consequently $\dH(\R2)\hookrightarrow \W(\R2)$. Moreover, this embedding is compact. If $\Phi$ belongs to a compact set $\Omega$, then $C(\Phi)$ can be replaced by $C(\Omega)$.
\end{lemma}
\begin{proof}
Fix $R>0$ and denote $B_R=B(0,R)$ the ball of radius $R$ centered at the origin. It is easy to see $\dH\to \W(B_R)$ is a compact embedding.

Fix any bounded sequence $\{f_k\}\subset \dH(\R2)$. Using a diagonal argument, we can extract a subsequence and a function $f\in \dH(\R2)$ such that for any $R>0$ it holds $f_k\to f$ in $\W(B_R)$-norm.  We want to prove that $f_k\to f$ in $\W(\R2)$-norm. 

 Suppose that $\chi:(0,\infty)\to (0,\infty)$ is a smooth non-decreasing function such that $\chi(r)=0$ if $r\in(0,2)$ and $\chi(r)=1$ is $r\geq 3$. 
Since \eqref{decay-4} implies that $|\nabla \Phi|^2(z)\leq C_\Phi (1+|z|)^{-4}$, then we can apply Hardy's inequality (see \cite[eq (4.19)]{adimurthi2010hardy}) to obtain
\begin{align*}
    \|f_k-f\|_{\W(B_R^c)}^2&\leq \int_{\R2}\chi(|z|)^2|f_k-f|^2|\nabla \Phi|^2\\
    &\leq CR^{-2}(\log R)^2\int_{\R2}\frac{\chi(|z|)^2|f_k-f|^2}{|z|^2(1+\log|z|)^2}\\
    &\leq CR^{-2}(\log R)^2\left[\|\nabla f_k-\nabla f\|_{L^2(\R2)}^2+\|f_k-f\|_{L^2(B_{3}\setminus B_{2})}^2\right].
\end{align*}
Given any $\varepsilon>0$, we choose $R$ sufficient large such that 
\begin{align*}
    \|f_k-f\|_{\W(B_R^c)}^2\leq \frac{\varepsilon}{2}+\|f_k-f\|_{L^2(B_{3}\setminus B_{2})}^2.
\end{align*}
Fixing such $R$, we have $\|f_k-f\|_{\W(B_R)}^2+\|f_k-f\|_{L^2(B_{3}\setminus B_{2})}^2<\varepsilon/2$ when $k$ is large enough. Combining this with the above inequality, one has 
\begin{align*}
    \|f_k-f\|_{\W}^2&=\|f_k-f\|_{\W(B_R)}^2+\|f_k-f\|_{\W(B_R^c)}^2<\varepsilon
\end{align*}
provided that $k$ is large enough. This completes the proof $f_k\to f$ in $\W(\R2)$ norm. Therefore the embedding is compact. The rest conclusion is easy to see.


\end{proof}

\begin{lemma}[Poincar\'e type inequality]\label{lem:poincare}
    Fix any harmonic map $\Phi$, $1\leq p< \infty$, there exists a constant $C_{\Phi,p}>0$ such that for any $v\in \dot H^1(\R2;\mathbb{R}^3)=\{v\in H^1_{loc}(\R2; \mathbb{R}^3):\int_{\R2}|\nabla v|^2<\infty\}$ with  $\int_{\R2}v|\nabla\Phi|^2=0$ one has 
    \begin{align} 
        \left(\int_{\R2}|v|^p|\nabla\Phi|^2\right)^{\frac{1}{p}}\leq C_{\Phi,p}\left(\int_{\R2}|\nabla v|^2\right)^{\frac{1}{2}}.
    \end{align}
    Moreover, if $\Phi\in \Omega$ a compact set of harmonic maps, then $C_{\Phi,p}$ can be replaced by some uniform constant $C_{\Omega,p}$. 
\end{lemma}
\begin{proof}
We can repeat the proof of Lemma \ref{lem:cpt} to see that $\dot H^1(\R2;\mathbb{R}^3)$ is compactly embedded in $\{v\in H^1_{loc}(\R2; \mathbb{R}^3):\int_{\R2}|v|^p|\nabla \Phi|^2<\infty\}$. Our conclusion follows from a standard contradiction argument. For instance, one can see \cite[section 5.8.1]{evans1998partial} and the argument in the following Lemma \ref{lem:eig-gap}.
\end{proof}

\begin{lemma}\label{lem:Lap-1}
The inverse operator $(|\nabla \Phi|^{-2}\Delta)^{-1}$ is a well-defined and continuous mapping from $\W(\R2)$ into $\H(\R2)$. 
\end{lemma}
\begin{proof}
Let $g\in \H(\R2)$ and $f\in \W(\R2)$. Applying H\"older's inequality, we obtain,
\begin{align} 
\begin{split}
    \langle f,g\rangle_{\W}=\int_{\R2}f\cdot g|\nabla \Phi|^2\leq  \|f\|_{\W}\| g\|_{\W}.
\end{split} 
\end{align}
As a consequence,
\[f\mapsto \langle f,\cdot\rangle_{\W}\in (\H)'\]
is continuous and injective. By the Riesz representation theorem, there exists a unique continuous linear map $T:\W(\R2)\to \H(\R2)$ such that for any $f\in \W$ and $g\in\H$
\begin{align} 
\int_{\R2}f\cdot g|\nabla \Phi|^2=\int_{\R2}\nabla T(f): \nabla g=-\int_{\R2}\Delta T(f)\cdot g.
\end{align}
Thus $-\Delta T(f)=f|\nabla \Phi|^2$, which implies $(|\nabla \Phi|^{-2}\Delta)^{-1}=-T$ satisfies the conclusion.
\end{proof}

\begin{lemma}\label{lem:eig-gap}
Suppose $\Omega$ is a compact set of degree $d$ harmonic maps. There exists $\mu=\mu(\Omega)>1$ such that for any $\Phi\in \Omega$
\begin{align}
    \int_{\R2}|\nabla v|^2\geq \mu\int_{\R2}|\nabla \Phi|^2|v|^2,\quad v\in (\ker \L[\Phi])^{\perp}\subset \H.
\end{align}
Here the orthogonality is with respect to $(\cdot,\cdot)_{\H}$.
\end{lemma}
\begin{proof}
Lemma \ref{lem:minimizer} says that each harmonic map achieves minimal Dirichlet energy in its homotopy class, therefore the second variation of $\mathcal{E}(u)$ is non-negative. It follows from \eqref{lin-2nd} that a density argument that 
\begin{align} 
    \int_{\R2}|\nabla v|^2-|\nabla \Phi|^2|v|^2\geq 0, \quad v\in\H(\R2).
\end{align}
The equality holds if and only if $v\in \ker \L[\Phi]$. Thus 
\begin{align} 
    \inf_{v\in \H}\frac{\int_{\R2}|\nabla v|^2}{\int_{\R2}|\nabla \Phi|^2|v|^2}=1.
\end{align}
Using Lemma \ref{lem:cpt} and Lemma \ref{lem:Lap-1}, $(|\nabla \Phi|^{-2}\Delta)^{-1}:\W\to \W$ is a compact self-adjoint operator, thus its spectrum is discrete. By the min-max characterization of eigenvalues, there exists $\mu_2(\Phi)>1$ such that
\begin{align} 
    \int_{\R2}|\nabla v|^2\geq \mu_2(\Phi)\int_{\R2}|\nabla \Phi|^2|v|^2,\quad v\in (\ker\L[\Phi])^\perp.
\end{align}
Here the orthogonality is with respect to $(\cdot,\cdot)_{\H}$.

Next, we want to show that $\exists\, \mu>1$ such that $\mu_{2}(\Phi)>\mu$ for all $\Phi\in \Omega$. Suppose not, then there exists a sequence of $\Phi_k=\mathcal{S}(p_k/q_k)\in \Omega$, $v_k\in (\ker\L[\Phi_k])^\perp$ such that 
\begin{align}\label{nab<1/k}
    \int_{\R2}|\nabla v_k|^2\leq (1+\frac{1}{k})\int_{\R2}|\nabla \Phi_k|^2|v_k|^2.
\end{align}
 Since $\Omega$ is a compact set, going to a subsequence if necessary, we can assume $p_k\to p_*$ and $q_k\to q_*$. After rescaling, we assume that $\|v_k\|_{\mathcal{H}_{\Phi_k}}=1$ and $\Phi_*=\mathcal{S}(p_*/q_*)$. Then \eqref{nab<1/k} implies $v_k\in \dot H^1(\R2)$ for any $k\geq 1$. Similar to Lemma \ref{lem:cpt}, $\dot H^1(\R2)$ is compactly embedded in $\mathcal{W}=\{v\in L_{loc}^1(\R2;\mathbb{R}^3):\int_{\R2}(1+|z|)^{-4}|v|^2<\infty\}$. Therefore $v_k$ all belong to the weighted space $\mathcal{H}$ defined in \eqref{def:calH}. Then there exists $v_*\in \mathcal{H}$ such that, going to a subsequence if necessary, $v_k\to v_*$ weakly in $\mathcal{H}$ and strongly in $\mathcal{W}$. Recall that \eqref{nablaS2} implies $ |\nabla \Phi|\leq C(1+|z|)^{-2}$ uniformly for all $\Phi\in \Omega$. Then $v_*\in \mathcal{H}_{\Phi_*}$. Since $p_k\to p_*$ and $q_k\to q_*$ then $\nabla \Phi_k\to \nabla \Phi_*$ a.e. Therefore by dominated convergence theorem,
 \begin{align} 
    \int_{\R2}|\nabla \Phi_*|^2|v_*|^2= \lim_{k\to \infty}\int_{\R2}|\nabla \Phi_k|^2|v_k|^2.
\end{align}
Taking limit in \eqref{nab<1/k}, we obtain from the above fact that
 \begin{align} 
     \int_{\R2}|\nabla v_*|^2=\int_{\R2}|\nabla\Phi_*|^2|v_*|^2=\frac12.
 \end{align}
 Thus $v_*\in \ker\L[\Phi_*]$ and $v_*\neq 0$. On the other hand, it follows from the non-degeneracy of $\L[\Phi_k]$ (see Lemma \ref{prop:non-deg}) that $\ker\L[\Phi_k]$ consists of smooth vector fields obtained from taking derivative of coefficients of $p_k$ and $q_k$. Moreover, Remark \ref{rem:K} concludes that these vector fields depend smoothly on $p_k,q_k$ and thus converge to a vector field in $\ker\L[\Phi_*]$ as $k\to \infty$. 
Since $v_k\in (\ker\L[\Phi_k])^{\perp}$ means
\begin{align} 
    \int_{\R2}\nabla v_k:\nabla K_k+|\nabla \Phi_k|^2v_k\cdot K_k=0,\quad K_k\in \ker\L[\Phi_k].
\end{align}
Letting $k\to \infty$, we obtain $(v_*,K_*)_{\mathcal{H}_{\Phi_*}}=0$ for some $K_*\in \ker \L[\Phi_*]$. Conversely for any $K_*$, we can choose a sequence of vector fields in $\ker \L[\Phi_k]$ converging to it. Therefore $v_*\in (\ker\L[\Phi_*])^{\perp}$. This contradicts the previous fact. The lemma is proved.
\end{proof}

The following Lemma is crucial for our local stability theorem. The proof here is a slight modification of that in \cite{bernand2021quantitative} for degree $\pm 1$ case.
\begin{lemma}\label{lem:p-Sob-Poin}
Suppose $\Omega$ is a compact set of degree $d$ harmonic maps and $p>1$. There exist two constants $\eta_4(\Omega)$ and $C_{\Omega,p}$ such that if $u\in \dot H^1(\R2;\S2)$ and $\|u-\Phi\|_{\dot H^1}\leq\eta_4(\Omega)$ for some harmonic map $\Phi\in \Omega$, then
\begin{align}\label{p-Sob-Poin}
    \int_{\R2}|u-\Phi|^p|\nabla\Phi|^2\leq C_{\Omega,p}\left(\int_{\R2}|\nabla (u-\Phi)|^2\right)^{\frac{p}{2}}.
\end{align}
\end{lemma}

\begin{proof}
We will prove the theorem assuming $u$ is smooth. By a density result of \citet{schoen1983boundary} it holds for any $u\in \dot H^1(\R2;\S2)$. Using \eqref{nab_u-deg}, we have $|\nabla u|^2\geq 2|u\cdot (u_x\times u_y)|$. Thus we obtain
\begin{align}
\begin{split}\label{u-nab_u}
    &\left|\int_{\R2}u\left(2|u\cdot(u_x\times u_y)|-|\nabla u|^2\right)\right|\leq \int_{\R2}\left(|\nabla u|^2-2|u\cdot(u_x\times u_y)|\right)\\
    \leq&\min\left\{\int_{\R2}\left(|\nabla u|^2+ 2u\cdot(u_x\times u_y)\right),\int_{\R2}\left(|\nabla u|^2- 2u\cdot(u_x\times u_y)\right)\right\}\\
    = &\ 2[\mathcal{E}(u)-4\pi|\deg(u)|].
    \end{split}
\end{align}
Take $\eta_4(\Omega)=\min\{\eta_1(\Omega),1/100\}$ (see \eqref{H1-contin}), then we have $\deg(u)=\deg(\Phi)$. Notice that 
\begin{align*}
    2[\mathcal{E}(u)-4\pi|\deg(u)|]=\int_{\R2}|\nabla u|^2-|\nabla\Phi|^2=\int_{\R2}|\nabla(u-\Phi)|^2+2\nabla (u-\Phi):\nabla\Phi.
\end{align*}
Since $|\nabla\Phi|\leq C(1+|z|)^{-2}$ as $z\to \infty$, we can apply integration by parts for the last term to see that
\begin{align}
2\int_{\R2}\nabla(u-\Phi):\nabla \Phi=2\int_{\R2}\Phi\cdot (u-\Phi)|\nabla\Phi|^2=-\int_{\R2}|u-\Phi|^2|\nabla \Phi|^2,
\end{align}
where we have used the fact $2\Phi\cdot(u-\Phi)=-|u-\Phi|^2$ in the last step. It follows that
\begin{align}\label{2Eu-degu}
    2[\mathcal{E}(u)-4\pi|\deg(u)|]=\int_{\R2}|\nabla (u-\Phi)|^2-|u-\Phi|^2|\nabla \Phi|^2.
\end{align}
Plugging this back to \eqref{u-nab_u}, we obtain
\begin{align}\label{u||-nab_u2}
    \left|\int_{\R2}u\left(2|u\cdot(u_x\times u_y)|-|\nabla u|^2\right)\right|\leq \int_{\R2}|\nabla(u-\Phi)|^2.
\end{align}
Since $u\cdot u_x=0$ and $u\cdot u_y=0$, the two vectors $u$ and $u_x\times u_y$ are parallel. Therefore, we have 
\begin{align} 
    |u\cdot (u_x\times u_y)|^2=|u_x\times u_y|^2=|u_x|^2|u_y|^2-(u_x\cdot u_y)^2.
\end{align}
Thus if we think of $u$ as a mapping from $\R2$ to $\mathbb{R}^3$, then $|u\cdot(u_x\times u_y)|$ is the modulus of the Jacobian of $u$. By the coarea formula (see Theorem 2.17 in \cite{ambrosio2000functions}), one has 
\begin{align}\label{u||-1}
    \int_{\R2}u|u\cdot (u_x\times u_y)|=\int_{\S2}z\mathcal{H}^0(\{u^{-1}(z)\})d\mu_{\S2}
\end{align}
where $\mathcal{H}^0$ is the 0-dimensional Hausdorff measure, i.e. counting the number of points. By Sard's theorem, for almost all $z\in \S2$, $z$ is a regular value of $u$. For such $z$, one has $\mathcal{H}^0(\{u^{-1}(z)\})\geq |\deg(u)|$ (see \cite[pg 95]{outerelo2009mapping}). By symmetry of the sphere we get
\begin{align}
\begin{split}\label{u||-2}
    \left|\int_{\S2}z\mathcal{H}^0(\{u^{-1}(z)\})d\mu_{\S2}\right|=&\left|\int_{\S2}z\left(\mathcal{H}^0(\{u^{-1}(z)\})-|\deg(u)|\right)d\mu_{\S2}\right|\\
    \leq &\int_{\S2}\left(\mathcal{H}^0(\{u^{-1}(z)\})-|\deg(u)|\right)d\mu_{\S2}.
    \end{split}
\end{align}
Using coarea formula again
\begin{align}
\begin{split}\label{u||-3}
    \int_{\S2}\left(\mathcal{H}^0(\{u^{-1}(z)\})-|\deg(u)|\right)d\mu_{\S2}=&\int_{\R2}|u\cdot (u_x\times u_y)|dx-4\pi|\deg(u)|\\
    \leq &\ \frac{1}{2}\int_{\R2}|\nabla u|^2-4\pi|\deg(u)|
    \end{split}
\end{align}
where we used \eqref{nab_u-deg}. Now one concatenates \eqref{u||-1}, \eqref{u||-2}, \eqref{u||-3} and \eqref{2Eu-degu} to get
\begin{align}\label{ineq-5-14}
   \left| \int_{\R2}u|u\cdot(u_x\times u_y)| \right|\leq \frac{1}{2}\int_{\R2}|\nabla(u-\Phi)|^2.
\end{align}
Using \eqref{ineq-5-14} and \eqref{u||-nab_u2}, we get
\begin{align} 
   \left| \int_{\R2}u|\nabla u|^2\right|\leq 2\int_{\R2}|\nabla(u-\Phi)|^2.
\end{align}
Applying $|\nabla u|^2-|\nabla \Phi|^2=2\nabla \Phi:\nabla(u-\Phi)+|\nabla(u-\Phi)|^2$ and Cauchy-Schwarz inequality, one gets
\begin{align} 
    \left|\int_{\R2}u(|\nabla u|^2-|\nabla \Phi|^2)\right|\leq 4\sqrt{2\pi}\left(\int_{\R2}|\nabla (u-\Phi)|^2\right)^{\frac12}+\int_{\R2}|\nabla(u-\Phi)|^2.
\end{align}
Since $\int_{\R2}|\nabla(u-\Phi)|^2<\eta_4(\Omega)\leq 1/100$, then 
\begin{align} 
    \left|\int_{\R2}u|\nabla\Phi|^2\right|\leq 8\sqrt{\pi}\left(\int_{\R2}|\nabla (u-\Phi)|^2\right)^{\frac12}.
\end{align}
Since $\int_{\R2}\Phi|\nabla\Phi|^2=\int_{\R2}\Delta\Phi=0$, then 
\begin{align}\label{u-Phi|nabPhi|}
    \left|\int_{\R2}(u-\Phi)|\nabla\Phi|^2\right|\leq 8\sqrt{\pi}\left(\int_{\R2}|\nabla (u-\Phi)|^2\right)^{\frac12}.
\end{align}
Now we apply Lemma \ref{lem:poincare} to get 
\begin{align} 
    \int_{\R2}|u-\Phi|^p|\nabla \Phi|^2\leq C_{\Omega,p}\left(\int_{\R2}|\nabla(u-\Phi)|^2\right)^{\frac{p}{2}}+C_{\Omega,p}\left|\int_{\R2}(u-\Phi)|\nabla\Phi|^2\right|^p.
\end{align}
Inserting \eqref{u-Phi|nabPhi|}, one gets \eqref{p-Sob-Poin}. 
\end{proof}

Now we can prove the local stability as claimed in the introduction. Our proof follows closely to the one  in \cite[ Lemma 4.4]{bernand2021quantitative}.
\begin{proof}[\bf Proof of Theorem \ref{thm:main-2}] Since the Dirichlet energy is invariant under M\"obius transformations, without loss of generality we assume $F=id_{\mathbb{R}^2}$. Otherwise one can work with $u\circ F^{-1}$. 

Since $\Omega$ is a compact, then  $\mathcal{A}_{\Omega}^m=\{(p,q):p\text{ is monic and }\mathcal{S}(p/q)\in \Omega\}$ is a compact subset of $\mathcal{A}_d$. Therefore there exists a $\varepsilon_1(\Omega)$ such that the following subset is also compact. 
\begin{align}
    (\mathcal{A}_\Omega^m)_{\varepsilon_1}=\{(\tilde p/\tilde q):\exists\, (p,q)\in \mathcal{A}_\Omega^m \text{ such that } |\tilde p-p|_\infty+|\tilde q-q|_\infty\leq \varepsilon_1 \}.
\end{align}
Denote $\Omega_{\varepsilon_1}=\mathcal{S}((\mathcal{A}_\Omega^m)_{\varepsilon_1})$. It is a compact set of degree $d$ harmonic maps. 
We shall take
\begin{align}\label{choice-eta}
    \tilde \eta(\Omega)=\min\{\eta_1(\Omega_{\varepsilon_1}),\eta_2(\Omega_{\varepsilon_1},\varepsilon_1),\eta_3(\Omega_{\varepsilon_1},\varepsilon_1),\eta_4(\Omega_{\varepsilon_1})\}.
\end{align}
It follows from Lemma \ref{lem:inf-ach} that the infimum can be achieved. Let us assume it is achieved at $\Phi$. Then Corollary \ref{cor:pq-near} and Lemma \ref{lem:pq-near} imply that $\Phi\in\Omega_{\varepsilon_1}$.  Denote $\delta=\int_{\R2}|\nabla(u-\Phi)|^2$. We decompose $\zeta=u-\Phi$ into three parts
\begin{align}
    \zeta_{\parallel}=(\zeta\cdot\Phi)\Phi,\quad  \zeta_{K}=proj_{\ker \L[\Phi]}(\zeta-\zeta_{\parallel}),\quad \zeta^*=\zeta-\zeta_{\parallel}-\zeta_K.
\end{align}
Here the projection is with respect to the inner product of $\W$. Consequently 
\begin{align}\label{orth-W}
    \int_{\R2}\zeta_K\cdot \zeta^*|\nabla \Phi|^2=0.
\end{align}
Since $\Delta \zeta_K+2(\nabla \Phi:\nabla \zeta_K)\Phi+|\nabla \Phi|^2\zeta_K=0$,  the above orthogonality is also equivalent to the orthogonality in $\dot H^1$.
\begin{align}\label{orth-H1}
    \int_{\R2}\nabla \zeta_K:\nabla \zeta^*=0.
\end{align}
\begin{claim}\label{claim:zeta2} Assume $\delta<\tilde \eta(\Omega)$. There exists $C=C(\Omega)$ such that 
\begin{align}
    \int_{\R2}|\zeta|^2|\nabla\Phi|^2\leq \int_{\R2}| \zeta^*|^2|\nabla \Phi|^2+C\delta^4.
\end{align}
\end{claim}
Indeed, the observation starts from the following identity 
\begin{align}\label{zeta2-whole}
    \int_{\R2}|\zeta|^2|\nabla \Phi|^2=\int_{\R2}(|\zeta_{\parallel}|^2+|\zeta_K|^2+|\zeta^*|^2)|\nabla \Phi|^2
\end{align}
which follows from $(\zeta_K+\zeta^*)\cdot \zeta_{\parallel}=0$ and \eqref{orth-W}.

First, using the smallness assumption $\delta<\tilde \eta(\Omega)\leq \eta_4(\Omega_{\varepsilon_1})$ and $\zeta_{\parallel}=-\frac12|\zeta|^2\Phi$, we can apply Lemma \ref{lem:p-Sob-Poin} to get
\begin{align}\label{u-Phi<delta4}
    \int_{\R2}|\zeta_{\parallel}|^2|\nabla \Phi|^2=\frac14\int_{\R2}|u-\Phi|^4|\nabla\Phi|^2\leq C\delta^4.
\end{align}

Second, since $\|u-\Phi\|_{\dot H^1}$ achieves the infimum at $\Phi$, one has $\int_{\R2}\nabla\zeta:\nabla\zeta_K=0$. Since $\zeta=\zeta_{\parallel}+\zeta_{K}+\zeta^*$ and \eqref{orth-H1}, this is equivalent to 
\begin{align}
    \int_{\R2}|\nabla \zeta_K|^2=-\int_{\R2}\nabla\zeta_K:\nabla\zeta_{\parallel}.
\end{align}
Since $\zeta_K$ is smooth and decay fast (see Remark \ref{rem:K}), we may do integration by parts for the right hand side and use the facts $\zeta_K\in \ker \L[\Phi]$ and $\zeta_K\cdot\zeta_{\parallel}=0$ to get
\begin{align}
    \int_{\R2}|\nabla\zeta_K|^2=-2\int_{\R2}(\zeta_{\parallel}\cdot \Phi)(\nabla \Phi:\nabla \zeta_K)=\int_{\R2}|u-\Phi|^2(\nabla \Phi:\nabla \zeta_K).
\end{align}
We apply H\"older's inequality to the above equation and \eqref{u-Phi<delta4} to get
\begin{align}
    \int_{\R2}|\nabla \zeta_K|^2\leq C\delta^4.
\end{align}
Now Lemma \ref{lem:cpt} implies that
\begin{align}\label{zetaK-delta4}
    \int_{\R2}|\zeta_K|^2|\nabla\Phi|^2\leq C \delta^4.
\end{align}
Inserting \eqref{u-Phi<delta4} and \eqref{zetaK-delta4} into \eqref{zeta2-whole}, we can prove the Claim \ref{claim:zeta2}.

\begin{claim}\label{claim:nabla_zeta}
Assume $\delta<\tilde\eta(\Omega)$. There exists $C=C(\Omega)$ such that 
\begin{align}
    \int_{\R2}|\nabla\zeta|^2=\int_{\R2}|\nabla \zeta_{\parallel}|^2+|\nabla \zeta_K|^2+|\nabla \zeta^*|^2+2\nabla \zeta_{\parallel}:\nabla(\zeta-\zeta_{\parallel})
\end{align}
and 
\begin{align}
\left|\int_{\R2}2\zeta_{\parallel}:\nabla(\zeta-\zeta_{\parallel})\right|\leq C\delta^3.
\end{align}
\end{claim}
Indeed, the first identity follows from \eqref{orth-H1}. From $\zeta_{\parallel}=-\frac12|\zeta|^2\Phi$, we have $\partial_\alpha \zeta_{\parallel}^l=-\frac12|\zeta|^2\partial_\alpha\Phi^l-\frac12\Phi^l\partial_{\alpha}|\zeta|^2$, where $\alpha\in \{x,y\}$ and $l\in \{1,2,3\}$. Using this identity, we have 
\begin{align}
    \begin{split}
\int_{\R2}2\nabla\zeta_{\parallel}:\nabla(\zeta-\zeta_{\parallel})=-\int_{\R2}|\zeta|^2\partial_\alpha\Phi^l\partial^\alpha(\zeta-\zeta_{\parallel})_l+\Phi^l\partial_\alpha|\zeta|^2\partial^\alpha(\zeta-\zeta_{\parallel})_l\\
    =\int_{\R2}(\zeta-\zeta_{\parallel})_l\partial_\alpha|\zeta|^2\partial^\alpha \Phi^l-|\zeta|^2\partial_\alpha\Phi^l\partial^\alpha(\zeta-\zeta_{\parallel})_l  
    \end{split}
\end{align} 
where we have used $\partial_{\alpha}[\Phi\cdot(\zeta-\zeta_{\parallel})]=0$. Since $|\nabla \Phi|\leq C(1+|z|)^{-2}$, we may integrate by parts for the first term to get
\begin{align}
    \int_{\R2}(\zeta-\zeta_{\parallel})_l\partial_\alpha|\zeta|^2\partial^\alpha \Phi^l=-\int_{\R2}|\zeta|^2\partial_\alpha\Phi^l\partial^\alpha(\zeta-\zeta_{\parallel})_l.
\end{align}
Therefore,
\begin{align}
     \int_{\R2}2\nabla\zeta_{\parallel}:\nabla(\zeta-\zeta_{\parallel})=-2\int_{\R2}|\zeta|^2\partial_\alpha\Phi^l\partial^\alpha(\zeta-\zeta_{\parallel})_l.
\end{align}
Since $\partial_\alpha \Phi_l\partial^\alpha \zeta_{\parallel}^l=-\frac{1}{2}\partial_\alpha \Phi_l\partial^\alpha[|u-\Phi|^2\Phi^l]=-\frac{1}{2}|u-\Phi|^2||\nabla \Phi|^2$, 
\begin{align}
    \int_{\R2}2\nabla\zeta_{\parallel}:\nabla(\zeta-\zeta_{\parallel})=&-2\int_{\R2}|\zeta|^2\nabla\Phi:\nabla\zeta+\int_{\R2}|u-\Phi|^4|\nabla\Phi|^2.
\end{align}
Using H\"older's inequality and Lemma \ref{lem:p-Sob-Poin}, we obtain
\begin{align}
\left|\int_{\R2}2\nabla\zeta_{\parallel}:\nabla(\zeta-\zeta_{\parallel})\right|\leq 2\delta\left(\int_{\R2}|\zeta|^4|\nabla\Phi|^2\right)^{\frac12}+C\delta^4\leq C\delta^3.
\end{align}
Thus Claim \ref{claim:nabla_zeta} is proved.

Now we can use these two claims to prove the theorem. Recall \eqref{2Eu-degu} and Lemma \ref{lem:eig-gap}.
\begin{align*}
    \begin{split}
         \mathcal{E}(u)-4\pi|\deg(u)|=&\ \int_{\R2}|\nabla \zeta|^2-|\zeta|^2|\nabla \Phi|^2\\
    \geq &\  \int_{\R2}|\nabla \zeta^*|^2-|\zeta^*|^2|\nabla\Phi|^2+\int_{\R2}|\nabla \zeta_{\parallel}|^2+|\nabla\zeta_K|^2-C\delta^3-C\delta^4\\
    \geq &\ (1-\mu^{-1})\int_{\R2}|\nabla\zeta^*|^2+\int_{\R2}|\nabla \zeta_{\parallel}|^2+|\nabla\zeta_K|^2-C\delta^3-C\delta^4\\
    \geq &\ (1-\mu^{-1})\int_{\R2}|\nabla \zeta|^2-C\delta^3-C\delta^4=(1-\mu^{-1})\delta^2-C\delta^3-C\delta^4.
    \end{split}
\end{align*}
Choosing $\eta(\Omega)=\min\{\tilde\eta(\Omega), \frac{1}{4C}(1-\mu^{-1})\}$ where $C$ is obtained from the above line. If $\delta<\eta(\Omega)$, then 
\begin{align}
    \mathcal{E}(u)-4\pi|\deg(u)|\geq \frac12(1-\mu^{-1})\delta^2=\frac12(1-\mu^{-1})\int_{\R2}|u-\Phi|^2.
\end{align}
The theorem is proved.
\end{proof}

\section{Counterexample with degree 2}\label{sec:counterexample}
 In this section, we shall construct some example to fulfill \eqref{stab-high} and thus Theorem \ref{thm:main-1} is established. The process starts with a particular degree 2 harmonic map $\mathcal{S}((z-r-\im r)(z+r+\im r))$. Here $r$ will be chosen large enough to satisfy various conditions. We shall introduce some notations first. 
Denote $\vec{\alpha}=(\alpha_i)\in\mathbb{R}^{10}$ and
\begin{align}\label{def:Psi[]}
   \Psi[{\vec{\alpha}}](z)= \frac{(\alpha_1+\im \alpha_{2}) z^2+(\alpha_3+\im \alpha_4)z+(\alpha_5+\im \alpha_6)}{1-(\alpha_7+\im \alpha_8)z -(\alpha_9+\im \alpha_{10})z^2}.
\end{align}
We also define $K_i=K_i[\vec{\alpha}](z)$, for $i=1,\cdots,10$, as
\begin{align}
\begin{split}\label{def-Ki}
    K_i[\vec{\alpha}](z)=&r^{\beta_i}{\partial_{\alpha_i} \mathcal{S}(\Psi[\vec{\alpha}])}(z) \\
    =&\frac{r^{\beta_i}}{(1+|\Psi|^2)^2}\left( 2(1+|\Psi|^2)\partial_{\alpha_i} \Psi- 2\Psi\partial_{\alpha_i}|\Psi|^2, 2\partial_{\alpha_i}|\Psi|^2\right)[\vec{\alpha}](z)
    \end{split}
\end{align}
where $\beta_i=0$ if $i\in \{1,2,5,6\}$, $\beta_i=-1$ if $i\in \{3,4,7,8\}$, $\beta_i=-2$ if $i\in \{9,10\}$. The reason why we divide some $K_i$ by $r$ or $r^2$ is to make sure that $\mathcal{J}_{ij}$ (see \eqref{def:Jij}) is bounded above by some constant. It will be clear in Lemma \ref{lem:Jr-pt2}.  
The non-degeneracy result implies $dim_{\mathbb{R}}\ker \L[\mathcal{S}(\Psi[\vec{\alpha}])]=10$ when $\Psi[\vec{\alpha}](z)$ is an irreducible rational function. Actually we can prove directly $\{K_i:i=1,\cdots,10\}$ are linearly independent, therefore $\ker \L[\mathcal{S}(\Psi[\vec{\alpha}])]=\text{Span}_{\mathbb{R}}\{K_1,\cdots,K_{10}\}$.

\begin{lemma}\label{lem:K_i-ind}
Suppose $\vec{\alpha}$ satisfies that $\Psi[\vec{\alpha}](z)$ is an irreducible rational function and $\alpha_1+\im \alpha_2\neq 0$. Then  $\{K_i:i=1,\cdots,10\}$ is linearly independent.
\end{lemma}
\begin{proof}
Suppose $\sum_{i=1}^{10}c_iK_i[\vec{\alpha}](z)=0$ for some constants $c_i\in \mathbb{R}$.
Using the expression of $K_i$ in \eqref{def-Ki}, it implies that there exists $c_i$ such that 
\begin{align} 
    \sum_{i=1}^{10}c_i\partial_{\alpha_i}|\Psi|^2[\vec{\alpha}](z)=0,\quad \sum_{i=1}^{10}c_i\partial_{\alpha_i}\Psi[\vec{\alpha}](z)=0.
\end{align}
Dividing the second equation by $\Psi$, we get 
\begin{align} 
    \sum_{i=1}^{10}c_i\frac{\partial_{\alpha_i}\Psi}{\Psi}[\vec{\alpha}](z)=0,\quad z\in\mathbb{C}\setminus \{\text{zeros and poles of }\Psi[\vec{\alpha}]\}.
\end{align}
That is
\begin{align*} 
\begin{split}
    &\frac{(c_1+\im c_2)z^2+(c_3+\im c_4)z+(c_5+\im c_6)}{(\alpha_1+\im \alpha_{2})z^2+(\alpha_3+\im \alpha_4)z+(\alpha_5+\im \alpha_6)} +\frac{(c_7+\im c_{8})z+(c_9+\im c_{10})z^2}{1-(\alpha_7+\im \alpha_8)z-(\alpha_9+\im \alpha_{10})z^2}
    =0.
    \end{split}
\end{align*}
Suppose $\{\zeta_1,\zeta_2\}$ are the roots of $(\alpha_1+\im \alpha_2)z^2+(\alpha_3+\im \alpha_4)+(\alpha_5+\im \alpha_6)=0$. Consider the roots of $1-(\alpha_7+\im \alpha_8)z-(\alpha_9+\im \alpha_{10})z^2=0$. It possibly has two roots, say $\{\zeta_3,\zeta_4\}$, or one root, or no root. In the first case,
since $\Psi[\vec{\alpha}](z)$ is irreducible, then  $\{\zeta_1,\zeta_2\}\cap \{\zeta_3,\zeta_4\}=\emptyset$. Let $z\to \zeta_3$ or $\zeta_4$, then we must have $c_7+\im c_8=c_9+\im c_{10}=0$. Consequently $c_1+\im c_{2}=c_3+\im c_4=c_5+\im c_6=0$. The other two cases are similar to prove.
\end{proof}
Denote $\mathcal{J}[\vec{\alpha}]=(\mathcal{J}_{ij})_{1\leq i,j\leq 10}$ where 
\begin{align}\label{def:Jij}
    \mathcal{J}_{ij}[\vec{\alpha}]=\int_{\R2}|\nabla \mathcal{S}(\Psi[\vec{\alpha}])|^2K_i[\vec{\alpha}]\cdot K_j[\vec{\alpha}].
\end{align}

\begin{lemma}\label{lem:J-ind}
Suppose $\vec{\alpha}$ satisfies that $\Psi[\vec{\alpha}](z)$ is an irreducible rational function and $\alpha_1+\im \alpha_2\neq 0$. Then $\mathcal{J}[\vec{\alpha}]$ is a positive definite matrix for any $r>0$.
\end{lemma}
\begin{proof} 
 Assume not, then there exists some $\vec{c}=(c_1,\cdots,c_{10})$ such that 
\begin{align} 
    0\geq \sum_{i,j=1}^{10}\mathcal{J}_{ij}[\vec{\alpha}]c_ic_j=\int_{\R2}|\nabla \mathcal{S}(\Psi[\vec{\alpha}])|^2|c_iK_i[\vec{\alpha}]|^2.
\end{align}
This implies $\sum_{i=1}^{10}c_iK_i=0$. However, this contradicts the linear independence of $K_i$.
\end{proof}

Since $\L[\mathcal{S}(\Psi[\vec{\alpha}])](K_i[\vec{\alpha}])=0$ and $K_i\cdot \mathcal{S}(\Psi[\vec{\alpha}])=0$, one has 
\begin{align}\label{J-2nd}
\int_{\R2}\nabla K_i[\vec{\alpha}]: \nabla K_j[\vec{\alpha}]=\int_{\R2}|\nabla \mathcal{S}(\Psi[\vec{\alpha}])|^2K_i[\vec{\alpha}]\cdot K_j[\vec{\alpha}]=\mathcal{J}_{ij}[\vec{\alpha}].
\end{align}

For any $r>0$, we denote $\vec{\alpha}_r=(1,0,0,0,0, 2r^2,0,0,0,0)$. Then   
\begin{align*}
    \Psi[\vec{\alpha}_r](z)=(z-r-\im r )(z+r+\im r).
\end{align*}
We will write $K_i^r=K_i[\vec{\alpha}_r]$ and $\J=\mathcal{J}[\vec{\alpha}_r]$ for short. In the following, one will see that $\J$ plays an important role in the analysis near the harmonic map $\mathcal{S}(\Psi[\vec{\alpha}_r])$. It is necessary to have a more detailed knowledge of each entry of $\J$, at least the leading orders of them as $r\to \infty$.

\begin{proposition}\label{prop:det-Jr}
 After some row and column permutations, we can represent $\J$ as a block diagonal matrix.
\begin{align}\label{J-switch}
    \J\sim \frac{16\pi}{3}diag\{A_1,A_2,A_3,A_4\}
\end{align}
with 
\begin{align*}
    A_1=&\begin{pmatrix}2&-4& \alpha_0\\-4&\lambda_3&\gamma_0\\\alpha_0 &\gamma_0& 4\end{pmatrix}+O(r^{-6}),\quad A_2=\begin{pmatrix}2& 4& -\alpha_0\\4& \lambda_3& \gamma_0\\ -\alpha_0&\gamma_0&4\end{pmatrix}+O(r^{-6}), \\
    A_3=&\begin{pmatrix}\lambda_1&\beta_0-r^{-2}\\\beta_0-r^{-2}&\lambda_2\end{pmatrix}+O(r^{-6}),\quad  A_4=\begin{pmatrix}\lambda_1&r^{-2}-\beta_0\\ r^{-2}-\beta_0&\lambda_2\end{pmatrix}+O(r^{-6}),
\end{align*}
where $A_1$ corresponds to $i,j\in \{1,10,6\}$, $A_2$ corresponds to $i,j\in \{2,9,5\}$, $A_3$ corresponds to $i,j\in \{3,8\}$, and $A_4$ corresponds to $i,j\in \{4,7\}$.
Here $\alpha_0\approx\beta_0\approx\gamma_0=O(r^{-\frac92})$, $\lambda_1=8+\frac14r^{-4}$, $\lambda_2=4+\frac{1}{2}r^{-4}$ and $\lambda_3=8+4r^{-4}$. 
\end{proposition}
\begin{remark}
One can compute the determinant of $\J$.
  \begin{align}\label{det-Jr}
    \begin{split}
        \det\J=&\ \frac{16^{10}\pi^{10}}{3^{10}}(\det A_1)(\det A_2)(\det A_3)(\det A_4)\\
        =&\ \frac{16^{10}\pi^{10}}{3^{10}}\left(32r^{-4}+O(r^{-6})\right)^2\left(32 +4r^{-4}+O(r^{-6})\right)^2\\
        =&\ \frac{2^{60}\pi^{10}}{3^{10}}r^{-8}+O(r^{-10}).
    \end{split}
\end{align}
One can see the degenerate tendency of $\J$ as $r\to \infty$.
\end{remark}

The proof of this lemma needs some involved integration. We need to expand the integrand to the third order to prove the result. Instead of diving into massive computations, we defer the proof of it to the next section and continue the main thread of our construction.  Let 
\begin{align}\label{def:calH}
    \mathcal{H}=\{u\in H_{loc}^1(\R2;\mathbb{R}^3):\int_{\R2}|\nabla u|^2+(1+|z|)^{-4}|u|^2<\infty\}.
\end{align} 
It is easy to see that $\mathcal{H}$ is a Hilbert space.


\begin{proposition}\label{prop:IFT}
Fix any $r>0$, there exists a $\varepsilon_2(r)$ and $\eta_5(r)$ such that for any $u:\R2\to \S2$ with $\|u-\mathcal{S}(\Psi[{\vec{\alpha}_r}])\|_{\mathcal{H}}<\eta_5(r)$, then there exists a unique $\vec{\alpha}=\vec{\alpha}(u)$ satisfying $\|\vec{\alpha}-\vec{\alpha}_r\|<\varepsilon_2(r)$ and 
\begin{align} 
    \int_{\R2} \nabla u: \nabla K=0,\quad \text{for } K\in \ker \L[{\mathcal{S}(\Psi[{\vec{\alpha}}])}].
\end{align}
\end{proposition}
\begin{proof}
Define the following map
\begin{align} 
    \begin{split}
       F:\mathcal{H}\times \mathbb{R}^{10}&\to \mathbb{R}^{10}\\
    (u,\vec{\alpha})&\mapsto\left(\int_{\R2} \nabla u:\nabla K_1[\vec{\alpha}],\cdots,\int_{\R2} \nabla u:\nabla K_{10}[\vec{\alpha}]\right). 
    \end{split}
\end{align} 
Such map $F$ is well-defined because $u\in \mathcal{H}$. It is easy to see $F$ is smooth on $\vec{\alpha}$ because $K_i$ depends on $\vec{\alpha}$ smoothly. $F$ is at least $C^1$ on $u$.

For any $r>0$, there exists $\varepsilon_2(r)$ such that if $\|\vec{\alpha}-\vec{\alpha}_r\|<\varepsilon_2(r)$, then $\alpha$ satisfies the assumption of Lemma \ref{lem:K_i-ind}. Therefore $\int_{\R2}|\nabla \mathcal{S}(\Psi[\vec{\alpha}])|^2=16\pi$. Differentiating on $\vec{\alpha}$, it infers that
\begin{align} 
    \int_{\R2}\nabla \mathcal{S}(\Psi[\vec{\alpha}]):\nabla K_i[\vec{\alpha}]=0,\quad i=1,\cdots,10.
\end{align}
Equivalently, this is $F(\mathcal{S}(\Psi[\vec{\alpha}]),\vec{\alpha})=0$ and $F(u,\vec{\alpha})=F(u-\mathcal{S}(\Psi[\vec{\alpha}]),\vec{\alpha})$. We intend to apply implicit function theorem to $F$ at $(\mathcal{S}(\Psi[\vec{\alpha}_r]),\vec{\alpha}_r)$. The Jacobian matrix with respect to $\vec{\alpha}$ at  $(\mathcal{S}(\Psi[\vec{\alpha}_r]),\vec{\alpha}_r)$ is 
\begin{align} 
  \frac{\partial F}{\partial \vec{\alpha}}( (\mathcal{S}(\Psi[\vec{\alpha}_r]),\vec{\alpha}_r))=\left(\int_{\R2}\nabla K_i^r:\nabla K_j^r\right)_{1\leq i,j\leq 10}= \mathcal{J}^r.
\end{align}
Here we have used \eqref{J-2nd}. Lemma \ref{lem:J-ind} says that such Jacobian is non-degenerate. Therefore, using the implicit function theorem, there exist $\eta_5(r)$ and $\varepsilon_2(r)$ small enough such that if $\|u-\mathcal{S}(\Psi[\vec{\alpha}_r])\|_{\mathcal{H}}<\eta_6(r)$, then there exists a unique $\vec{\alpha}=\vec{\alpha}(u)$ such that $|\vec{\alpha}-\vec{\alpha}_r|<\varepsilon_2(r)$ and $F(u,\vec{\alpha})=0$ .
\end{proof}



Introduce the cut-off function
\begin{align} 
    \Theta^r(z)=\begin{cases}1 & |z|<{r}^{\frac12},\\2-2\log(|z|)/\log r& r^{\frac12}\leq |z|\leq r,\\0& |z|>r. \end{cases}
\end{align}
Define
\begin{align}\label{def:fr}
    f^r(x,y)=  \Theta^r(x+\im y-r-\im r)-\Theta^r(x+\im y+r+\im r)
\end{align}
and
\begin{align}\label{def:KK}
      \KK=  2K_2^r-K_9^r.
\end{align}
See the explicit formulae of $K_2^r,K_9^r$ in section \ref{sec:cal-Jac}. 
\begin{lemma}\label{lem:nablafK}
For $r>0$ large enough, we have 
\begin{align}
    \int_{\R2}|\nabla (f^r\KK)|^2=\frac{64\pi}{3}r^{-4}+O(|\log r|^{-1}r^{-4}), \\
    \int_{\R2}|\nabla \mathcal{S}(\Psi[\vec{\alpha}_r])|^2|f^r\KK|^2=\frac{64\pi}{3}r^{-4}+O(r^{-6}). 
\end{align}
\end{lemma}

 Denote
\begin{align}\label{def:pj}
    p_j=\int_{\R2}\nabla (f^r\KK)
    :\nabla K_j^rdxdy.
\end{align}
\begin{lemma}\label{lem:pj}
We have
\begin{align*}  
\begin{split}
    &p_1=-2p_{10}=O(r^{-\frac{13}{2}}),\quad p_2=O(r^{-6}),\quad p_3=-p_4 =O(r^{-\frac92}),\quad p_5=O(r^{-\frac92}),\\
    &p_6=O(r^{-\frac92}),\quad p_7=p_8=-\frac{16\pi}{3}r^{-4}+O(r^{-6}) ,\quad p_9=O(r^{-6}).
    \end{split}
\end{align*}
Consider the solution $\vec{c}=(c_1,\cdots,c_{10})^T$ of $\J\vec{c}=\vec{p}$ where $\vec{p}=(p_1,\cdots,p_{10})^T$,
\begin{align}\label{c_j}
  \begin{split}
       & c_1=2c_{10}+O(r^{-\frac{13}{2}})=O(r^{-\frac{5}{2}}), \quad c_5=O(r^{-\frac92}),\\ &c_2=-2c_9+O(r^{-6})=O(r^{-2}),\quad c_6=O(r^{-\frac92}),\\
       & c_3=-c_4=O(r^{-\frac92}),\quad c_7=c_8=-\frac{1}{4}r^{-4}+O(r^{-6}).
    \end{split}
\end{align} 
\end{lemma}
Lemma \ref{lem:nablafK} and Lemma \ref{lem:pj} are crucial for our construction. Since it requires technical computations, we also postpone the proofs of them to the next section.

\begin{proposition}\label{prop:orth-h}
For any $r>0$, there exists $\varepsilon_3(r)$ with the following significance. For any $\varepsilon<\varepsilon_3(r)$, there exists $\{h_1,\cdots,h_{10}\}$ which are $\dot H^1(\R2)\cap L^\infty(\R2)$ functions and depend on $\varepsilon,r$ continuously such that $u=\varepsilon {h}^i {K}_i^r+\sqrt{1-\varepsilon^2|h^iK_i^r|^2}\mathcal{S}(\Psi[\vec{\alpha}_r])$ satisfies 
\begin{align}\label{u:K-orth}
    \int_{\R2}\nabla u:\nabla K_i^r=0,\quad i=1,\cdots,10.
\end{align}
Furthermore, 
\begin{align} 
    &\int_{\R2}|\nabla(h^iK_i^r)|^2=\frac{64\pi}{3}r^{-4}+O(|\log r|^{-1}r^{-4})+O(\varepsilon),\label{hK-norm}\\
    &\int_{\R2}|\nabla(h^iK_i^r)|^2-\int_{\R2}|\nabla \mathcal{S}(\Psi[\vec{\alpha}_r])|^2|h^iK_i^r|^2=O(|\log r|^{-1}r^{-4})+O(\varepsilon).\label{hK-norm-2}
\end{align}
\end{proposition}
\begin{proof}
We can take $\vec{h}=(h_1,\cdots,h_{10})$ where
\begin{align} \label{def:hiKi}
    h^iK_i^r=f^r\KK-c^iK_i^r
\end{align}
with $c_i$ to be determined. Here we use the Einstein summation convention for $i$. 
Define a map
\begin{align} 
    F:\mathbb{R}_+\times \mathbb{R}^{10}&\to \mathbb{R}^{10}\\
    (\varepsilon,\vec{c}\,)&\to \left(\int_{\R2} \nabla v:\nabla K_1^r,\cdots, \int_{\R2}\nabla v:\nabla K_{10}^r\right)
\end{align}
where  
\begin{align} 
    v= h^iK_i^r-\frac{\varepsilon |h^iK_i^r|^2}{\sqrt{1-|\varepsilon h^iK_i^r|^2}+1}\mathcal{S}(\Psi[{\vec{\alpha}_r}]).
\end{align}
The map $F$ is well-defined because $K_i^r$ and $\mathcal{S}(\Psi[\vec{\alpha}_r])$ both belong to $\dot{H}^1(\mathbb{R}^2)$. For $\varepsilon$ and $|\vec{c}|$ small, $F$ is a smooth map. 

At $\varepsilon=0$, $F(0,\vec{c})=0$ if and only if 
\begin{align} 
    \J\vec{c}=\vec{p}
\end{align}
where $\vec{p}=(p_1,\cdots,p_{10})$, where $p_{j}$ is defined in \eqref{def:pj}. Since $\J$ is non-degenerate, using Lemma \ref{lem:pj}, the above equation has a unique solution \eqref{c_j}. We denote it as $\vec{c}_*=(c_1^*,\cdots,c_{10}^*)$. The Jacobian of $F$ at $(0,\vec{c}_*)$ with respect to $\varepsilon$ is 
\begin{align} 
    (\partial_{c_i}F^j)(0,\vec{c}_*)=- \J .
\end{align}
By the implicit function theorem, there exists $\varepsilon_3(r)$ such that for any $0\leq \varepsilon<\varepsilon_3(r)$ there exists $\vec{c}=\vec{c}(\varepsilon)=\vec{c}_*+O(\varepsilon)$ satisfies $F(\varepsilon,\vec{c})=0$. Since \eqref{def:hiKi} and $f^r\in \dot H^1(\R2)\cap L^\infty(\R2)$, then $h_i\in \dot{H}^1(\R2)\cap L^\infty(\R2)$ for $1\leq i\leq 10$. Using the form of $v$, one readily check $u=\varepsilon {h}^i {K}_i^r+\sqrt{1-\varepsilon^2|h^iK_i^r|^2}\mathcal{S}(\Psi[\vec{\alpha}_r])$ satisfies \eqref{u:K-orth}. 

To establish \eqref{hK-norm}, we compute explicitly
\begin{align}\label{inthK2}
\begin{split}
    \int_{\R2}|\nabla (h^iK_i^r)|^2
    =&\int_{\R2}|\nabla (f^r\KK)|^2-2\sum_{i=1}^{10}c_i\int_{\R2}\nabla (f^r\KK):\nabla K_i^r+\sum_{i,j=1}^{10}c_{i}c_j\J_{ij}\\
    =&\ \frac{64\pi}{3}r^{-4}-2\sum_{i=1}^{10}c_{i}^*p_i+\sum_{i,j=1}^{10}c_{i}^*c_{j}^*\J_{ij}+O(|\log r|^{-1}r^{-4})+O(\varepsilon)
    \end{split}
\end{align}
where we have used Lemma \ref{lem:nablafK} and $\vec{c}=\vec{c}_*+O(\varepsilon)$. 
By Lemma \ref{lem:pj}, we have
\begin{align} \label{cipi-1}
   2\sum_{i=1}^{10}c_i^*p_i =O(r^{-8}).
\end{align}
To compute $\sum c_i^*c_j^*\J_{ij}$, we shall use Proposition \ref{prop:det-Jr} to see that $\J$ can written as a block diagonal matrix. Combining \eqref{c_j}, Proposition \ref{prop:det-Jr}, we have 
\begin{align}
\begin{split}
    \sum_{i,j\in\{1,6,10\}}c_i^*c_j^*\J_{ij}=&\ (c_1^*)^2\J_{11}+2c_1^*c_{10}^*\J_{1,10}+(c_{10}^*)^2\J_{10,10}+O(r^{-8})\\
    =&\ \frac{32\pi}{3}\left((c_1^*)^2-2(c_1^*)^2+(c_1^*)^2\right)+O(r^{-8})=O(r^{-8})
    \end{split}
\end{align}
where we have used $c_1^*=2c_{10}^*+O(r^{-\frac{13}{2}})=O(r^{-\frac{5}{2}})$ in the second line. Similarly one can use  \eqref{c_j}, Proposition \ref{prop:det-Jr}  and $c_2^*=-2c_9^*+O(r^{-6})=O(r^{-2})$ to derive
\begin{align}\label{cicj-3}
    \sum_{i,j\in\{2,9,5\}}c_i^*c_j^*\J_{ij}=O(r^{-8}), \quad \sum_{i,j\in\{3,8\}}c_i^*c_j^*\J_{ij}=\sum_{i,j\in\{4,7\}}c_i^*c_j^*\J_{ij}=O(r^{-8}).
\end{align}
Plugging the equations \eqref{cipi-1}-\eqref{cicj-3} to \eqref{inthK2}, we can get \eqref{hK-norm}.

To establish \eqref{hK-norm-2}, we shall use $ h^iK_i^r=f^r\KK-c^iK_i^r$ and $K_i^r\in \ker \L[\mathcal{S}(\Psi[\vec{\alpha}_r])]$ to derive 
\begin{align*} 
    \int_{\R2}|\nabla(h^iK_i^r)|^2-\int_{\R2}|\nabla \mathcal{S}(\Psi[\vec{\alpha}_r])|^2|h^iK_i^r|^2= \int_{\R2}|\nabla (f^r\KK)|^2-|\nabla \mathcal{S}(\Psi[\vec{\alpha}_r])|^2|f^r\KK|^2.
\end{align*}
Then \eqref{hK-norm-2} just follows from Lemma \ref{lem:nablafK}. 
\end{proof}

\begin{proposition}\label{prop:inf-h}
For any $r>0$, there exists a $\varepsilon_4(r)$ such that for $\varepsilon<\varepsilon_4(r)$ there exists $\{h_1,\cdots,h_{10}\}$ which are $\dot H^1(\R2)\cap L^\infty(\R2)$ functions and depend on $\varepsilon,r$ continuously such that for $u=\varepsilon {h}^i {K}_i^r+\sqrt{1-\varepsilon^2|h^iK_i^r|^2}\mathcal{S}(\Psi[\vec{\alpha}_r])$ the following infimum is achieved at $\mathcal{S}({\Psi[\vec{\alpha}_r]})$. 
\begin{align}\label{infimum}
    \inf_{(p,q)\in \mathcal{A}_2}\|u- \mathcal{S}(p/q)\|_{\dot H^1}^2=\| u- \mathcal{S}(\Psi[\vec{\alpha}_r])\|_{\dot H^1}^2=\varepsilon^2\|h^iK_i^r\|_{\dot H^1}^2+O_r(\varepsilon^3).
\end{align}
Here $|O_r(\varepsilon^3)|\leq C(r)\varepsilon^3$ as $\varepsilon\to 0$.
\end{proposition}
\begin{proof}
First we choose $\varepsilon_4(r)<\varepsilon_3(r)$. Then it follows from Proposition \ref{prop:orth-h} that one can find $\{h_1,\cdots,h_{10}\}$ which are $\dot H^1(\R2)\cap L^\infty(\R2)$ functions and depend on $\varepsilon,r$ continuously such that $u=\varepsilon {h}^i {K}_i^r+\sqrt{1-\varepsilon^2|h^iK_i^r|^2}\mathcal{S}(\Psi[\vec{\alpha}_r])$ satisfies 
\begin{align}\label{u:K-orth-2}
    \int_{\R2}\nabla u:\nabla K_i^r=0,\quad i=1,\cdots,10.
\end{align}
Taking $\varepsilon_4(r)$ even smaller, we can make  \begin{align}
    \|u-\mathcal{S}(\Psi[\vec{\alpha}_r])\|_{\mathcal{H}}<\min\{\eta_2(r),\eta_6(r),\eta_4(r,\varepsilon_2(r)/100)\}.
\end{align}
It follows from Lemma \ref{lem:inf-ach} that the infimum \eqref{infimum} can be achieved , say $\tilde\Phi=\mathcal{S}(\tilde p/\tilde q)$ for some canonical $(\tilde p,\tilde q)\in \mathcal{A}_2$. Corollary \ref{cor:pq-near} implies that all minimizers are near $\mathcal{S}(\Psi[\vec{\alpha}_r])$. More precisely
\begin{align}
    |\tilde p-(z-r-\im r)(z+r+\im r)|_\infty+|\tilde q-1|_\infty\leq\frac{\varepsilon_2(r)}{100}. 
\end{align}  
Thus we can assume $\tilde \Phi=\mathcal{S}(\Psi[\vec{\alpha}])$ for some $\vec{\alpha}$ satisfies $\|\vec{\alpha}-\vec{\alpha}_r\|<\varepsilon_2(r)$.

Since the infimum achieved at $\mathcal{S}(\Psi[\vec{\alpha}])$, then one has the orthogonality condition
\begin{align} 
    \int_{\R2}\nabla u:\nabla K=0,\quad\forall \, K\in \ker \L[\mathcal{S}(\Psi[\vec{\alpha}])].
\end{align}
However, Proposition \ref{prop:IFT} says that such  $\vec{\alpha}$ is unique if $\|\vec{\alpha}-\vec{\alpha}_r\|<\varepsilon_2(r)$. Our choice of $\vec{h}$ from Proposition \ref{prop:orth-h} makes sure that $\vec{\alpha}$ has to be $\vec{\alpha}_r$.

Finally, since $h^iK_i^r$ is smooth and bounded on $\mathbb{R}^2$, we can compute explicitly 
\begin{align} 
    \|u-\mathcal{S}(\Psi[\vec{\alpha}_r])\|_{\dot H^1}^2=\|\varepsilon  h^iK_i^r+O(\varepsilon^2|\vec{h}|^2)\|_{\dot H^1}^2=\varepsilon^2\|h^iK_i^r\|_{\dot H^1}^2+O_r(\varepsilon^3).
\end{align}
\end{proof}

Finally we can prove the main theorem of this section. 
\begin{proof}[\bf Proof of Theorem \ref{thm:main-1}] We shall take $u$ from Proposition \ref{prop:inf-h}. 
Since $h^iK_i\in \dot H^1(\R2;\mathbb{R}^3)\cap L^\infty(\R2;\mathbb{R}^3)$, we shall apply Lemma \ref{lem:2nd-vari} to get
\begin{align} \label{RHS-log}
\begin{split}
    \mathcal{E}(u)-4\pi|\deg (u)|&=\frac12\varepsilon^2\int_{\R2}|\nabla (h^iK_i^r)|^2-|\nabla \mathcal{S}(\Psi[\vec{\alpha}_r])|^2|h^iK_i^r|^2+O_r(\varepsilon^3)\\
    &= \varepsilon^2 O(|\log r|^{-1}r^{-4})+O_r(\varepsilon^3)
    \end{split}
\end{align}
where we have used \eqref{hK-norm-2} in the last step. 

On the other hand, it follows from Proposition \ref{prop:inf-h} and \eqref{hK-norm} that 
\begin{align} \label{LHS-4}
    \inf_{(p,q)\in \mathcal{A}_2}\|u- \mathcal{S}(p/q)\|_{\dot H^1}^2=\frac{64\pi}{3}r^{-4}\varepsilon^2+O_r(\varepsilon^3)+O(\varepsilon^2r^{-6}).
\end{align}
Now for any $M>0$, we choose $r$ large such that $64\pi/(3M)> O(|\log r|^{-1})$. \textbf{Fixing} such $r$, we can choose $\varepsilon$ small such that 
\begin{align} 
    \inf_{\vec{\alpha}\in \mathbb{R}^{10}}\|u- \mathcal{S}(\Psi[\vec{\alpha}])\|_{\dot H^1}^2> M(\mathcal{E}(u)-4\pi |\deg(u)|).
\end{align}
Thus the theorem is proved.
\end{proof}


\section{Explicit computation of the Jacobian}\label{sec:cal-Jac}

In this section, we shall compute the Jacobian $\J$ explicitly and give   proofs of Proposition \ref{prop:det-Jr}, Lemma \ref{lem:nablafK}, and Lemma \ref{lem:pj}. 

Now consider $\Psi[\vec{\alpha}]$ defined in \eqref{def:Psi[]}. Identifying $\mathbb{C}^5=\{(a,b,c,d,e)\}$ and $\mathbb{R}^{10}=\{\vec{\alpha}\}$ in the way of $a=\alpha_1+\im \alpha_2$, $b=\alpha_3+\im \alpha_4$, etc. In this notation, we rewrite it as 
\begin{align*}
   \Psi[{\vec{\alpha}}](z)= \frac{az^2+bz+c}{1-dz -ez^2}.
\end{align*}
Also we rewrite $\vec{\alpha}_r=(1,0,-(r+\im r)^2,0,0)\in \mathbb{C}^5$,  and denote $\Psi[{\vec{\alpha}_r}]=\psi_r$ for short  
\begin{align*}
    \psi_r(z)=(z-r-\im r )(z+r+\im r)=x^2-y^2+\im(2xy-2r^2).
\end{align*}
We always assume that  $r$ is large enough. We want to compute the vector fields $K_i$ defined in \eqref{def-Ki} at $\vec{\alpha}_r$ explicitly. One needs to differentiate $\mathcal{S}(\Psi[\vec{\alpha}])$ with respect to real  and imaginary part of $a,b,c,d,e$ at $ \vec{\alpha}_r$. 
Therefore we need to compute $\partial \Psi[\vec{\alpha}]$ and $\partial|\Psi[\vec{\alpha}]|^2$ at $\vec{\alpha}_r$. It is easy to see 
\begin{align*}
    \begin{split}
        {\partial_{a_1} \Psi[{\vec{\alpha}_r}]}&=\psi_r,\quad
        {\partial_{a_2} \Psi[{\vec{\alpha}_r}]}=\im \psi_r,\quad 
        {\partial_{b_1} \Psi[{\vec{\alpha}_r}]}=z,\quad
        {\partial_{b_2} \Psi[{\vec{\alpha}_r}]}=\im z,\quad {\partial_{c_1} \Psi[{\vec{\alpha}_r}]}=1,\\
         {\partial_{c_2} \Psi[{\vec{\alpha}_r}]}&=\im ,\ 
         {\partial_{d_1} \Psi[{\vec{\alpha}_r}]}=z\psi_r,\ 
         {\partial_{d_2} \Psi[{\vec{\alpha}_r}]}=\im z\psi_r,\ 
         {\partial_{e_1} \Psi[{\vec{\alpha}_r}]}=z^2\psi_r,\ 
         {\partial_{e_2} \Psi[{\vec{\alpha}_r}]}=\im z^2\psi_r.
    \end{split}
\end{align*}
Note that $\partial|\psi|^2=\partial (\psi\bar \psi)=2\Re (\psi\partial\bar \psi)=2\Re(\psi\overline{\partial\psi})$.  Then
\begin{align*}
        {\partial_{a_1} |\Psi[{\vec{\alpha}_r}]|^2}&=2|\psi_r|^2,\quad
         {\partial_{a_2} |\Psi[{\vec{\alpha}_r}]|^2}=0,\\
        {\partial_{b_1} |\Psi[{\vec{\alpha}_r}]|^2}&=2(x^3-2r^2y+xy^2),\quad
        {\partial_{b_2} |\Psi[{\vec{\alpha}_r}]|^2}=2 (y^3-2r^2x+x^2y),\\
         {\partial_{c_1} |\Psi[{\vec{\alpha}_r}]|^2}&=2x^2-2y^2,\quad
         {\partial_{c_2} |\Psi[{\vec{\alpha}_r}]|^2}=4xy-4r^2,\\
        {\partial_{d_1} |\Psi[{\vec{\alpha}_r}]|^2}&=2x|\psi_r|^2,\quad
         {\partial_{d_2} |\Psi[{\vec{\alpha}_r}]|^2}=-2y|\psi_r|^2,\\
         {\partial_{e_1} |\Psi[{\vec{\alpha}_r}]|^2}&=2(x^2-y^2)|\psi_r|^2,\quad
         {\partial_{e_2} |\Psi[{\vec{\alpha}_r}]|^2}=-4xy|\psi_r|^2.
\end{align*}
Then we have
\begin{align*}
        &K_1^r:={\partial_{a_1} u}=\zeta\left( (1-|\psi_r|^2)\psi_r, 2|\psi_r|^2\right),\\
         &K_2^r:={\partial_{a_2} u}=\zeta\left( (1+|\psi_r|^2)\im \psi_r, 0\right),\\
          &K_3^r:=r^{-1}{\partial_{b_1} u}=\zeta r^{-1}\left((1+|\psi_r|^2) z -4(x^3-2r^2y+xy^2)\psi_r, 4(x^3-2r^2y+xy^2)\right),\\
          &K_4^r:=r^{-1}{\partial_{b_2} u}=\zeta r^{-1}\left((1+|\psi_r|^2)\im z-4(y^3-2r^2x+x^2y)\psi_r, 4(y^3-2r^2x+x^2y)\right),\\
          &K_5^r:= {\partial_{c_1} u}=\zeta\left((1+|\psi_r|^2) -4(x^2-y^2)\psi_r, 4(x^2-y^2)\right),\\
          &K_6^r:={\partial_{c_2} u}=\zeta\left((1+|\psi_r|^2)\im-8(xy-r^2)\psi_r, 8(xy-r^2)\right),\\
          &K_7^r:=r^{-1}{\partial_{d_1} u}=\zeta r^{-1}\left((1+|\psi_r|^2)z\psi_r-2x|\psi_r|^2 \psi_r, 2x|\psi_r |^2\right),\\
           &K_8^r:=r^{-1}{\partial_{d_2} u}=\zeta r^{-1}\left((1+|\psi_r|^2)\im z\psi_r+2y|\psi_r|^2 \psi_r, -2y|\psi_r |^2\right),\\
            &K_9^r:=r^{-2}{\partial_{e_1} u}=\zeta r^{-2}\left((1+|\psi_r|^2)z^2\psi_r-2(x^2-y^2)|\psi_r|^2 \psi_r, 2(x^2-y^2)|\psi_r |^2\right),\\
          &K_{10}^r:=r^{-2}{\partial_{e_2} u}=\zeta r^{-2}\left((1+|\psi_r|^2)\im z^2\psi_r+4xy|\psi_r|^2 \psi_r, -4xy|\psi_r |^2\right).
\end{align*}
where $\zeta=2(1+|\psi_r|^2)^{-2}$.
We will introduce the notation $\I_{ij}=\frac14(1+|\psi_r|^2)^2K^r_i\cdot K^r_j$ to denote the inner product of $K^r_i$ and $K^r_j$.
We have
\begin{align*}
    &\I_{11}={|\psi_r|^2}=\I_{22},\quad \I_{12}=0,\quad\I_{13}={r^{-1} (x^3 - 2 r^2 y + x y^2)},\\
    &\I_{14}={r^{-1} (y^3 - 2 r^2 x + x^2 y)}, \quad  \I_{15}= { (x^2 - y^2)},\quad
    \I_{16}= {2 (x y - r^2)},\quad \I_{17}={r^{-1}x|\psi_r|^2},
    \\
    &\I_{18}=-{r^{-1}y|\psi_r|^2 },\quad  \I_{19}={r^{-2}(x^2-y^2)|\psi_r|^2},\quad
    \I_{1,10}= -{2r^{-2}xy |\psi_r|^2},\\
          & \I_{23}=r^{-1}[{2 r^2 x -  x^2y - y^3)}],\quad \I_{24}=r^{-1}[{ (x^3 - 2 r^2 y + x y^2)}],\quad \I_{25}= {2 (r^2 - x y)},\\
          & \I_{26}= { (x^2 - y^2)},\quad\I_{27}={r^{-1}y|\psi_r|^2},\quad \I_{28}={r^{-1}x|\psi_r|^2},\\
        &\I_{29}={2r^{-2}xy|\psi_r|^2},\quad \I_{2,10}={r^{-2}(x^2-y^2)|\psi_r|^2}.
    \end{align*}
\begin{align*}
        &\I_{33}={r^{-2}(x^2+y^2)}=\I_{44},\quad  \I_{34}=0,\quad \I_{35}={r^{-1}x},\\ &\I_{36}={r^{-1}y},\quad\I_{37}={r^{-2} (x^4 - y^4) },\quad\I_{38}={2 r^{-2}(r^2 - x y) (x^2 + y^2)},\\
        &\I_{39}={r^{-3} (x^2 + y^2) (x^3 + 2 r^2 y - 3 x y^2)},\quad \I_{3,10}=-{r^{-3} (x^2 + y^2) (2 r^2 x - 3 x^2 y + y^3)},\\
    & \I_{45}=-{r^{-1}y},\quad \I_{46}={r^{-1}x},\quad  \I_{47}=-{r^{-2} (r^2 - x y) (x^2 + y^2)} ,\quad  \I_{48}={r^{-2} (x^4 - y^4)} ,\\
         &\I_{49}=-{r^{-3} (x^2 + y^2) (2 r^2 x - 3 x^2 y + y^3)},\quad \I_{4,10}={ r^{-3}(x^2 + y^2) (x^3 + 2 r^2 y - 3 x y^2)} .
\end{align*}
\begin{align*}
    &\I_{55}=1=\I_{66},\ \I_{56}=0,\ \I_{57}={r^{-1} (x^3 + 2 r^2 y - 3 x y^2) },\  \I_{58}={r^{-1} (2 r^2 x - 3 x^2 y + y^3)} ,\\
    &\I_{59}={r^{-2} (x^4 + 4 r^2 x y - 6 x^2 y^2 + y^4)} ,\  \I_{5,10}={2r^{-2} (r^2 - 2 x y) (x^2 - y^2) },\\
   &\I_{67}=\- {r^{-1} (2 r^2 x - 3 x^2 y + y^3)} , \quad
         \I_{68}=\ {r^{-1} (x^3 + 2 r^2 y - 3 x y^2)} ,\\
         &\I_{69}=\ {2r^{-2} (r^2 - 2 x y) (y^2 - x^2)} ,\quad 
         \I_{6,10}=\ {r^{-2} (x^4 + 4 r^2 x y - 6 x^2 y^2 + y^4)} .
\end{align*}
\begin{align*}
    &\I_{77}={r^{-2}(x^2+y^2)|\psi_r|^2} ,\quad  \I_{78}=0,\\
    &\I_{79}={r^{-3}x(x^2+y^2)|\psi_r|^2} ,\quad\I_{7,10}={-r^{-3}y(x^2+y^2)|\psi_r|^2},\\
    &\I_{88}={r^{-2}(x^2+y^2)|\psi_r|^2} ,\quad\I_{89}={r^{-3}y(x^2+y^2)|\psi_r|^2} ,\quad\I_{8,10}={r^{-3}x(x^2+y^2)|\psi_r|^2},\\
    &\I_{99}={r^{-4}(x^2+y^2)^2|\psi_r|^2} ,\quad \I_{9,10}=0,\quad\I_{10,10}={r^{-4}(x^2+y^2)^2|\psi_r|^2} .
\end{align*}

The calculations here and later are a little bit tedious but still manageable with bare hands. To make the life easier, we perform them with the help of Mathematica\footnote{\url{https://www.wolframcloud.com/obj/bingomat/Published/deg2-HM-final.nb} \\A link for source code.}.

Using \eqref{pu}, we get
\begin{align}
    |\nabla \mathcal{S}(\Psi[\vec{\alpha}_r])|^2=\frac{4|\partial_x\psi_r|^2+4|\partial_y\psi_r|^2}{(1+|\psi_r|^2)^2}=\frac{32(x^2+y^2)}{(1+|\psi_r|^2)^2}.
\end{align}
Thus \eqref{def:Jij} implies that
\begin{align}\label{J-I}
    \J_{ij}=\int_{\R2}\frac{128(x^2+y^2)}{(1+|\psi_r|^2)^4}\I_{ij}dxdy.
\end{align}

\begin{lemma}\label{lem:Jr-pt1}
The following entries of $\J$ are all equal to $0$.
\begin{align*}
    &\J_{12},\J_{13},\J_{14}, \J_{15},\J_{17},\J_{18},\J_{19},\J_{23},\J_{24},\J_{26}, \J_{27},\J_{28},\J_{2,10},\\
    &\J_{34},\J_{35},\J_{36},\J_{37},\J_{39},\J_{3,10},\J_{45},\J_{46}, \J_{48}, \J_{49}, \J_{4,10},\\
    &\J_{56},\J_{57},\J_{58},\J_{5,10},
     \J_{67},\J_{68},\J_{69},\J_{78},\J_{79},\J_{7,10}, \J_{89},\J_{8,10},\J_{9,10}.
\end{align*}
\end{lemma}
\begin{proof}
These facts follow from the symmetries.  The integrand of  $\J_{13}$ is odd with respect to the operation $(x,y)\to (-x,-y)$  while the integration domain is even. Thus it is equal to 0. The same symmetry holds for
\begin{align*}
    &\J_{14},\J_{17},\J_{18},\J_{23},\J_{24},\J_{27},\J_{28},\J_{35},\J_{36},\J_{39},\J_{3,10},\\
    &\J_{45},\J_{46},\J_{49},\J_{4,10},\J_{57},\J_{58},\J_{67},\J_{68},\J_{79},\J_{7,10},\J_{89},\J_{8,10}.
\end{align*} 
On the other hand, if we switch $x$ and $y$, then we get $\J_{15}=-\J_{15}$. The same symmetry holds for
\begin{align*}
    \J_{19},\J_{26},\J_{2,10},\J_{37},\J_{48},\J_{5,10},\J_{69}.
\end{align*} 
It is easy to see $\J_{12}=\J_{34}=\J_{56}=\J_{78}=\J_{9,10}=0$ since $\I_{12}=\I_{34}=\I_{56}=\I_{78}=\I_{9,10}=0$.
\end{proof}

To compute other terms of $\J_{ij}$, we need the following lemma. 
\begin{lemma}\label{lem:pPsi}
Suppose $p(x,y,r)$ is a homogeneous polynomial on $x,y,r$ with degree $k\geq 0$. Assume $l\geq 3 $  and $l>\frac{k}{4}+\frac12$. For $r$ large, one has 
\begin{align*}
    \int_{\R2}\frac{p(x,y,r)}{(1+|\psi_r|^2)^l}=&\frac{\pi}{8(l-1)} [p(\mathbf{1})+p(-\mathbf{1})]r^{k-2}+\frac{\pi[\Delta_{x,y}p(\mathbf{1})
    +\Delta_{x,y}p(-\mathbf{1})]}{256(l-1)(l-2)}r^{k-6}\\&+\frac{\pi[- p_x(\mathbf{1})-p_y(\mathbf{1})+ p_x(-\mathbf{1})+p_y(-\mathbf{1})+p(\mathbf{1})+p(-\mathbf{1})]}{128(l-1)(l-2)}r^{k-6}\\
    &+O(r^{k-8})
\end{align*}
where $\mathbf{1}=(1,1,1)$ and $-\mathbf{1}=(-1,-1,1)$.
\end{lemma}

\begin{proof}
Let
\begin{align*}
    \Omega_1&=\{(x,y):(x-r)^2+(y-r)^2<r\},\\
    \Omega_2&=\{(x,y):(x+r)^2+(y+r)^2<r\},\\
    \Omega_3&=\{(x,y): (x-r)^2+(y-r)^2\geq r,|y|+|x|\leq 3r,x+y\geq 0\},\\
    \Omega_4&=\{(x,y): (x+r)^2+(y+r)^2\geq r,|y|+|x|\leq 3r,x+y\leq 0\},
\end{align*}
and $\Omega^{out}=\R2-\Omega_1-\Omega_2-\Omega_3-\Omega_4\subset \{(x,y):|x|\geq r,|y|\geq r\}$.  We shall reserve the notations of $\Omega_i$ for the rest of paper. See Figure \ref{fig:int-domain} for illustration. 
Denote 
\[Int(\Omega_i)=\int_{\Omega_i}\frac{p(x,y,r)}{(1+|\psi_r|^2)^l}dxdy.\]

\begin{figure}
    \centering
    \includegraphics[width=0.35\textwidth]{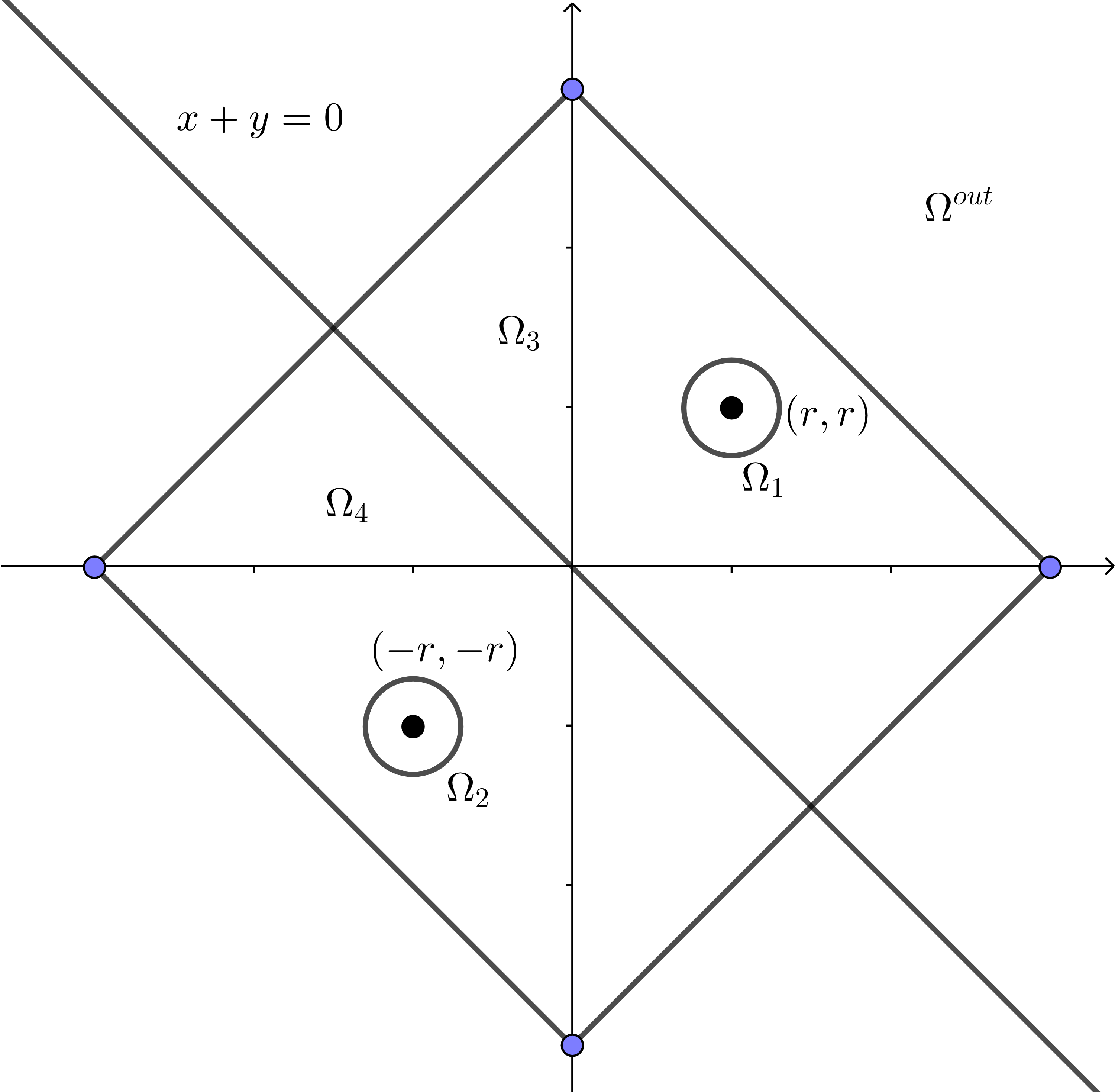}
    \caption{Integration domains $\Omega_i$.}
    \label{fig:int-domain}
\end{figure}

On $\Omega_1$, making a change of variable $x=r+r^{-1}\tilde x$, $y=r+r^{-1}\tilde y$, we have $\Omega_1=\{(\tilde x,\tilde y):\tilde x^2+\tilde y^2< r^3\}$. Since $|\psi_r|^2=[(x-r)^2+(y-r)^2][(x+r)^2+(y+r)^2]$, we rewrite it as
\begin{align} 
\begin{split}
    |\psi_r|^2=&(\tilde x^2+\tilde y^2)[(2+r^{-2}\tilde x)^2+(2+r^{-2}\tilde y)^2]\\
    =&(\tilde x^2+\tilde y^2)[8+4r^{-2}(\tilde x+\tilde y)+r^{-4}(\tilde x^2+\tilde y^2)].
    \end{split}
\end{align}
Denote $\tilde x^2+\tilde y^2=\rho^2$ and  $A=1+8(\tilde x^2+\tilde y^2)=1+8\rho^2$.  Since $r^{-2}|\tilde x|\leq r^{-1/2}$ and $r^{-2}|\tilde y|\leq r^{-1/2}$ in $\Omega_1$, then $4r^{-2}(\tilde x+\tilde y)+r^{-4}(\tilde x^2+\tilde y^2)\ll 8$. Then we can make the following expansion.
\begin{align}
\begin{split}\label{1+psi-exp1}
    &(1+|\psi_r|^2)^{-l}=\big(A+4r^{-2}\rho^2(\tilde x+\tilde y)+r^{-4}\rho^4\big)^{-l}\\
    =&\ A^{-l}-lA^{-l-1}[4r^{-2}\rho^2(\tilde x+\tilde y)+r^{-4}\rho^4]+\frac{l(l+1)}{2}A^{-l-2}[4r^{-2}\rho^2(\tilde x+\tilde y)+r^{-4}\rho^4]^2\\
    &+O(A^{-l-3}r^{-6}\rho^6(|\tilde x|^3+|\tilde y|^3)).
    \end{split}
\end{align}
Since $p$ is a homogeneous polynomial and $r^{-2}|\tilde x|\leq r^{-1/2}$ and $r^{-2}|\tilde y|\leq r^{-1/2}$, then
\begin{align}\label{p-exp1}
\begin{split}
    &p(x,y,r)= r^kp(1+r^{-2}\tilde x,1+r^{-2}\tilde y,1)\\
    =&\ r^kp(\mathbf{1})+r^{k-2}[p_x(\mathbf{1})\tilde x+p_y(\mathbf{1})\tilde y]+\frac12r^{k-4}[p_{xx}(\mathbf{1})\tilde x^2+p_{xy}(\mathbf{1})\tilde x\tilde y+p_{yy}(\mathbf{1})\tilde y^2]\\&+O(r^{k-6}(|\tilde x|^3+|\tilde y|^3)).
\end{split}
\end{align}
Here we write $(1,1,1)$ as $\mathbf{1}$ for short. Now we combine \eqref{1+psi-exp1} and \eqref{p-exp1} to get 
\begin{align} 
    \frac{p(x,y,r)}{(1+|\psi_r|^2)^l}=r^{k}\frac{p(\mathbf{1})}{[1+8\rho^2]^l}+r^{k-2}B_1+r^{k-4}B_2+O(r^{k-6}A^{-l}(|\tilde x|^3+|\tilde y|^3)),
\end{align}
where 
\begin{align} 
    B_1=\frac{p_x(\mathbf{1})\tilde x+p_y(\mathbf{1})\tilde y}{[1+8\rho^2]^l}-\frac{4p(\mathbf{1})l\rho^2(\tilde x+\tilde y)}{[1+8\rho^2]^{l+1}}
\end{align}
and 
\begin{align} 
\begin{split}
    B_2=&\ \frac{[p_{xx}(\mathbf{1})\tilde x^2+p_{xy}(\mathbf{1})\tilde x\tilde y+p_{yy}(\mathbf{1})\tilde y^2]}{2[1+8\rho^2]^{l}}-\frac{4\rho^2l(\tilde x+\tilde y)[p_x(\mathbf{1})\tilde x+p_y(\mathbf{1})\tilde y]}{[1+8\rho^2]^{l+1}}\\
    &-\frac{p(\mathbf{1})l\rho^4}{[1+8\rho^2]^{l+1}}+\frac{8p(\mathbf{1}){l(l+1)}\rho^4(\tilde x+\tilde y)^2}{[1+8\rho^2]^{l+2}}.
    \end{split}
\end{align}
Making a change of variable $\tilde x=\rho\cos \theta$, $\tilde y=\rho \sin\theta$, we have 
\begin{align*}
     Int(\Omega_1)
    =&\ 2\pi r^{k-2}\int_{0}^{r^{\frac{3}{2}}}\frac{p(\mathbf{1})\rho}{[1+8\rho^2]^l}d\rho +r^{k-6}\int_{\Omega_1}B_2 d\tilde xd\tilde y +O(r^{k-8})\int_0^{r^{\frac{3}{2}}} \frac{\rho^3}{[1+8\rho^2]^l}d\rho\\
    =&\ 2\pi p(\mathbf{1})r^{k-2}\int_{0}^{r^{\frac{3}{2}}}\frac{\rho}{(1+8\rho^2)^l}d\rho+r^{k-6}\int_{\Omega_1}B_2d\tilde xd\tilde y+O(r^{k-8}).
\end{align*}
Here we have used the assumption $l\geq  3$ and $\int_{\Omega_1}B_1d\tilde xd\tilde y=0$. One can compute
\begin{align*}
    \int_{\Omega_1}B_2d\tilde xd\tilde y=&\ \frac{\pi}{2}\Delta_{x,y} p(\mathbf{1})\int_{0}^{r^{\frac{3}{2}}}\frac{\rho^3}{[1+8\rho^2]^l}d\rho-4\pi[p_x(\mathbf{1})+p_y(\mathbf{1})]l\int_{0}^{r^{\frac{3}{2}}}\frac{\rho^5}{[1+8\rho^2]^{l+1}}d\rho\\
    &\ -2\pi p(\mathbf{1})l\int_0^{r^{\frac{3}{2}}}\frac{\rho^5}{[1+8\rho^2]^{l+1}}d\rho+16\pi p(\mathbf{1})l(l+1)\int_0^{r^{\frac{3}{2}}}\frac{\rho^7}{[1+8\rho^2]^{l+2}}d\rho.
\end{align*}
Now some elementary integration shows 
\begin{align*} 
\begin{split}
 \int_0^{r^{\frac{3}{2}}} \frac{\rho}{[1+8\rho^2]^l}d\rho=&\ \frac{1}{16(l-1)}+O(r^{-6}),\\
    \int_{0}^{r^{\frac{3}{2}}}\frac{\rho^3}{[1+8\rho^2]^l}d\rho=&\ \frac{1}{128(l-1)(l-2)}+O(r^{-3}),\\
    \int_{0}^{r^{\frac{3}{2}}}\frac{\rho^5}{[1+8\rho^2]^{l+1}}d\rho=&\ \frac{1}{512l(l-1)(l-2)}+O(r^{-3}),\\
    \int_0^{r^{\frac{3}{2}}}\frac{\rho^7}{[1+8\rho^2]^{l+2}}d\rho=&\ \frac{3}{4096(l+1)l(l-1)(l-2)}+O(r^{-3}).   
\end{split}
\end{align*}
Plugging in these identities to $\int_{\Omega_1}B_2d\tilde xd\tilde y$ to get
\begin{align*}
    &\int_{\Omega_1}B_2d\tilde xd\tilde y\\
    =&\frac{\pi\Delta_{x,y}p(\mathbf{1})}{256(l-1)(l-2)}-\frac{2\pi[p_x(\mathbf{1})+p_y(\mathbf{1})]+\pi p(\mathbf{1})}{256(l-1)(l-2)}+\frac{3\pi p(\mathbf{1})}{256(l-1)(l-2)}+O(r^{-3})\\
    =&\frac{\pi[\Delta_{x,y}p(\mathbf{1})-2 p_x(\mathbf{1})-2p_y(\mathbf{1})+2p(\mathbf{1})]}{256(l-1)(l-2)}+O(r^{-3}).
\end{align*}
Therefore
\begin{align}\label{IOmega1}
    Int(\Omega_1)=\frac{\pi p(\mathbf{1})}{8(l-1)}r^{k-2}+\frac{\pi[\Delta_{x,y}p(\mathbf{1})-2 p_x(\mathbf{1})-2p_y(\mathbf{1})+2p(\mathbf{1})]}{256(l-1)(l-2)}r^{k-6}+O(r^{k-8}).
\end{align}
By symmetry, we have corresponding equality for $Int(\Omega_2)$,
\begin{align}\label{IOmega2}
Int(\Omega_2)=&\frac{\pi p(-\mathbf{1})}{8(l-1)}r^{k-2}+\frac{\pi[\Delta_{x,y}p(-\mathbf{1})+2 p_x(-\mathbf{1})+2p_y(-\mathbf{1})+2p(-\mathbf{1})]}{256(l-1)(l-2)}r^{k-6}+O(r^{k-8})
\end{align}
where $p(-\mathbf{1})=(-1,-1,1)$.

Now in $\Omega_3$, we have $r^{-2}|\tilde x|+r^{-2}|\tilde y|\leq C $ and $|\psi_r|^2\geq C^{-1} r^2[(x-r)^2+(y-r)^2]$ for some uniform constant $C$. Let $x=r+r^{-1}\tilde x$ and $y=r+r^{-1}\tilde y$, one gets
\begin{align}
\begin{split}\label{IOmega3}
|Int(\Omega_3)|\leq&\  Cr^k\int_{\Omega_3}|\psi_r|^{-2l}dxdy\leq Cr^{k-2}\int_{r^{\frac{3}{2}}}^\infty [\tilde x^2+\tilde y^2]^{-l}d\tilde x d\tilde y\\
\leq&\  Cr^{k-2}\int_{r^\frac{3}{2}}^{+\infty}\rho^{-2l+1}d\rho\leq Cr^{k-2} r^{-3l+3}=O(r^{k-8}).
    \end{split}
\end{align}
The same estimate holds for $Int(\Omega_4)$. 
 In $\Omega^{out}$, we make a change of variable $x=r\tilde x$, $y=r \tilde y$,
\begin{align}
\begin{split}\label{IOmegaout}
    |Int(\Omega^{out})|&\leq \int_{\Omega^{out}}|p(x,y,r)||\psi_r|^{-2l}dxdy\\
    &\leq \int_{|\tilde x|\geq 1,|\tilde y|\geq 1}\frac{|p(\tilde x,\tilde y,1)|r^{k-4l+2}}{[((\tilde x-1)^2+(\tilde y-1)^2)((\tilde x-1)^2+(\tilde y-1)^2)]^{l}}d\tilde xd\tilde y\\
    &\leq r^{k-4l+2}\int_{|\tilde x|\geq 1,|\tilde y|\geq 1}\frac{|\tilde x|^k+|\tilde y|^k+C}{[\tilde x^2+\tilde y^2]^{2l}}d\tilde xd\tilde y\leq Cr^{k-4l+2}=O(r^{k-10})
    \end{split}
\end{align}
provided $4l>k+2$.

Collecting the results of $Int(\Omega_i)$, $i=1,2,3,4$ and $Int(\Omega^{out})$, we get 
\begin{align} 
\begin{split}\label{rk-8}
   Int(\R2)=&\frac{\pi}{8(l-1)} [p(\mathbf{1})+p(-\mathbf{1})]r^{k-2}+\frac{\pi[\Delta_{x,y}p(\mathbf{1})
    +\Delta_{x,y}p(-\mathbf{1})]}{256(l-1)(l-2)}r^{k-6}\\&+\frac{\pi[- p_x(\mathbf{1})-p_y(\mathbf{1})+ p_x(-\mathbf{1})+p_y(-\mathbf{1})+p(\mathbf{1})+p(-\mathbf{1})]}{128(l-1)(l-2)}r^{k-6}\\
    &+O(r^{k-8}). 
\end{split} 
\end{align}
 This proves our conclusion.
\end{proof}
\begin{corollary}\label{cor:int-psi}
Suppose $p(x,y,r)$ is a homogeneous polynomial on $x,y,r$ with degree $k\geq 0$. Assume $l\geq 4 $  and $l>\frac{k}{4}+\frac32$. For $r$ large, one has 
\begin{align} 
\begin{split}
    &\int_{\R2}\frac{p(x,y,r)|\psi_r|^2}{(1+|\psi_r|^2)^l}\\
    =&\frac{\pi}{8(l-1)(l-2)} [p(\mathbf{1})+p(-\mathbf{1})]r^{k-2}+\frac{\pi[\Delta_{x,y}p(\mathbf{1})
    +\Delta_{x,y}p(-\mathbf{1})]}{128(l-1)(l-2)(l-3)}r^{k-6}\\&+\frac{\pi[- p_x(\mathbf{1})-p_y(\mathbf{1})+ p_x(-\mathbf{1})+p_y(-\mathbf{1})+p(\mathbf{1})+p(-\mathbf{1})]}{64(l-1)(l-2)(l-3)}r^{k-6}+O(r^{k-8}). 
\end{split}
\end{align}
\end{corollary}
\begin{proof}
Notice the following equality
\begin{align} 
    \int_{\R2}\frac{p(x,y,r)|\psi_r|^2}{[1+|\psi_r|^2]^l}=\int_{\R2}\frac{p(x,y,r)}{[1+|\psi_r|^2]^{l-1}}-\int_{\R2}\frac{p(x,y,r)}{[1+|\psi_r|^2]^l},
\end{align}
one can apply Lemma \ref{lem:pPsi}.
\end{proof}
 Next we shall prove that if $p(1,1,1)=p(-1,-1,1)=0$ then the remainder of term in \eqref{rk-8} could be improved.
\begin{lemma}\label{lem:pPsi-2}
Suppose that $p(x,y,r)$ is a homogeneous polynomial on $x,y,r$ with degree $k\geq 0$. Assume $p(1,1,1)=p(-1,-1,1)=0$,  $l\geq 3$, $l>\frac{k}{4}+\frac12$, then
\begin{align*}
\begin{split}
    &\int_{\R2}\frac{p(x,y,r)}{(1+|\psi_r|^2)^l}dxdy\\
    =&\frac{\pi[\Delta_{x,y} p(\mathbf{1})+\Delta_{x,y}p(-\mathbf{1})-2p_x(\mathbf{1})-2p_y(\mathbf{1})+2p_x(-\mathbf{1})+2p_y(-\mathbf{1})] r^{k-6}}{256(l-1)(l-2)}+O(r^{k-\frac{17}{2}}).
    \end{split}
\end{align*}
\end{lemma}
\begin{proof}
 We shall use the notations in the proof of Lemma \ref{lem:pPsi} and refine the proof there. In $\Omega_1$, 
 \begin{align}
\begin{split}\label{1+psi-exp2}
    &(1+|\psi_r|^2)^{-l}=\big(A+4r^{-2}\rho^2(\tilde x+\tilde y)+r^{-4}\rho^4\big)^{-l}\\
    =&A^{-l}-lA^{-l-1}[4r^{-2}\rho^2(\tilde x+\tilde y)+r^{-4}\rho^4]+\frac{l(l+1)}{2}A^{-l-2}[4r^{-2}\rho^2(\tilde x+\tilde y)+r^{-4}\rho^4]^2\\
    &-\frac{l(l+1)(l+2)}{6}A^{-l-3}[4r^{-2}\rho^2(\tilde x+\tilde y)+r^{-4}\rho^4]^3+O(A^{-l-4}r^{-8}\rho^8(|\tilde x|^4+|\tilde y|^4)|).
    \end{split}
\end{align}
Using our assumption, one has 
\begin{align}\label{p-exp2}
\begin{split}
    &p(x,y,r)=r^kp(1+r^{-2}\tilde x,1+r^{-2}\tilde y,1)\\
    =&r^{k-2}[p_x(\mathbf{1})\tilde x+p_y(\mathbf{1})\tilde y]+\frac12r^{k-4}[p_{xx}(\mathbf{1})\tilde x^2+p_{xy}(\mathbf{1})\tilde x\tilde y+p_{yy}(\mathbf{1})\tilde y^2]\\&+\frac16 r^{k-6}Q+O(r^{k-8}(|\tilde x|^4+|\tilde y|^4)) 
\end{split}
\end{align}
where 
\begin{align} 
    Q=p_{xxx}(\mathbf{1})\tilde x^3+3p_{xxy}(\mathbf{1})\tilde x^2\tilde y+3p_{xyy}(\mathbf{1})\tilde x\tilde y^2+p_{yyy}(\mathbf{1})\tilde y^3.
\end{align}
Now we combine \eqref{1+psi-exp2} and \eqref{p-exp2} to get 
\begin{align} 
    \frac{p(x,y,r)}{(1+|\psi_r|^2)^l}=r^{k-2}B_1+r^{k-4}B_2+r^{k-6}B_3+O(r^{k-8}A^{-l}(|\tilde x|^4+|\tilde y|^4))
\end{align}
where 
\begin{align*}\begin{split}
    B_1=&\ \frac{p_x(\mathbf{1})\tilde x+p_y(\mathbf{1})\tilde y}{[1+8\rho^2]^l},\\
    B_2=&\ \frac{[p_{xx}(\mathbf{1})\tilde x^2+p_{xy}(\mathbf{1})\tilde x\tilde y+p_{yy}(\mathbf{1})\tilde y^2]}{2[1+8\rho^2]^{l}}-\frac{4\rho^2l(\tilde x+\tilde y)[p_x(\mathbf{1})\tilde x+p_y(\mathbf{1})\tilde y]}{[1+8\rho^2]^{l+1}},\\
    B_3=&\ \frac{Q}{6[1+8\rho^2]^l}-\frac{l\rho^4[p_{x}(\mathbf{1})\tilde x+p_y(\mathbf{1})\tilde y]+2l\rho^2(\tilde x+\tilde y)[p_{xx}(\mathbf{1})\tilde x^2+p_{xy}(\mathbf{1})\tilde x\tilde y+p_{yy}(\mathbf{1})\tilde y^2]}{[1+8\rho^2]^{l+1}}.
    \end{split}
\end{align*}
Using $\int_{\Omega_1}B_3d\tilde x d\tilde y=0=\int_{\Omega_1}B_1d\tilde xd\tilde y$, we can get the following estimate
\begin{align}\label{IOmega1-2nd}
    Int(\Omega_1)=\frac{\pi[\Delta_{x,y}p(\mathbf{1})-2 p_x(\mathbf{1})-2p_y(\mathbf{1})]}{256(l-1)(l-2)}r^{k-6}+O(r^{k-10}).
\end{align}
Similar equality hold for $Int(\Omega_2)$. In $\Omega_3$, using our assumption, 
\begin{align*}
\begin{split}
|Int(\Omega_3)|\leq&\  Cr^{k-2}\int_{\Omega_3}(|\tilde x|+|\tilde y|)|\psi_r|^{-2l}dxdy\leq Cr^{k-4}\int_{r^{\frac{3}{2}}}^\infty [\tilde x+\tilde y][\tilde x^2+\tilde y^2]^{-l}d\tilde x d\tilde y\\
\leq&\  Cr^{k-4}\int_{r^\frac{3}{2}}^{+\infty}\rho^{-2l+2}d\rho\leq Cr^{k-4} r^{-3l+\frac92}=O(r^{k-\frac{17}{2}}).
    \end{split}
\end{align*}
The same estimate holds for $Int(\Omega_4)$. For $\Omega^{out}$, \eqref{IOmegaout} still holds.
Collecting the results of $Int(\Omega_i)$, $i=1,2,3,4$ and $Int(\Omega^{out})$, we get the conclusion.
\end{proof}

\begin{lemma}\label{lem:Jr-pt2}
We have 
\begin{align*}
    &\J_{11}=\J_{22}=\frac{32\pi}{3}+O(r^{-6}), \quad\J_{1,10}=-\frac{64\pi}{3} +O(r^{-6}),\\
    & \J_{16}=-\J_{25}=O(r^{-\frac92}) ,\quad \J_{29}=\frac{64\pi}{3}+O(r^{-6}),\\
    &\J_{33}=\J_{44}=\frac{128}{3}\pi+\frac{4\pi}{3}r^{-4}
    +O(r^{-6}),  \quad \J_{38}=-\J_{47}=-\frac{16\pi}{3}r^{-2}+O(r^{-\frac92}), \\
  &\J_{55}=\J_{66}=\frac{64\pi}{3}+O(r^{-6}),\quad \J_{59}= \J_{6,10}=O(r^{-\frac92}),\\
     &\J_{77}=\J_{88}=\frac{64\pi}{3}+\frac{8\pi}{3}r^{-4}+O(r^{-6}), \\ &\J_{99}=\J_{10,10}=\frac{128\pi}{3}+\frac{64\pi}{3}r^{-4}+O(r^{-6}) .
\end{align*}
\end{lemma}
\begin{proof}
Recall \eqref{J-I}.
Applying Corollary \ref{cor:int-psi} with $k=2, l=4$ and $p(x,y,r)=128(x^2+y^2)$, we have
\begin{align} 
    \J_{11}=\int_{\R2}\frac{128(x^2+y^2)|\psi_r|^2}{(1+|\psi_r|^2)^4}dxdy=\frac{32\pi}{3}+O(r^{-6}).
\end{align}
In the same way, one can compute $\J_{1,10},\J_{22},\J_{29},\J_{33},\J_{44},\J_{55},\J_{66},\J_{77},\J_{88},\J_{99},\J_{10,10}$.

 Applying Lemma \ref{lem:pPsi-2} with $k=4, l=4$ and $p(x,y,r)=256(x^2+y^2)(xy-r^2)$, we have
 \begin{align} 
    \J_{16}=\int_{\R2}\frac{256(x^2+y^2)(xy-r^2)}{(1+|\psi_r|^2)^4}dxdy= O(r^{-\frac92}).
\end{align}
In the same way, one can compute $\J_{25},\J_{38},\J_{47},\J_{59},\J_{6,10}$.
\end{proof}

Next, we give proofs of Proposition \ref{prop:det-Jr} and two key lemmas required in section \ref{sec:counterexample}.
\begin{proof}[\bf Proof of Proposition \ref{prop:det-Jr}]
Combining Lemma \ref{lem:Jr-pt1} and Lemma \ref{lem:Jr-pt2}, we know that $\J$ has the specific form. 
\end{proof}
\begin{proof}[\bf Proof of Lemma \ref{lem:pj}]
It is easy to see that 
\begin{align}
     \int_{\R2}\nabla (f^rK_i^r):\nabla K_j^rdxdy=\int_{\R2}|\nabla \mathcal{S}(\Psi[\vec{\alpha}_r])|^2f^rK_i^r\cdot K_j^rdxdy,\quad  1\leq i,j\leq 10. 
\end{align}
We shall adopt the notation
\begin{align}
    Int_{ij}(\Omega)=\int_{\Omega}|\nabla \mathcal{S}(\Psi[\vec{\alpha}_r])|^2f^r(x,y)K_i^r\cdot K_j^rdxdy. 
\end{align}
Recall the definition of $p_j$ in \eqref{def:pj}. Then 
\begin{align*}
    p_j=2\left(\sum_{k=1}^4Int_{2j}\left(\Omega_k\right)+Int_{2j}\left(\Omega^{out}\right)\right)-\left(\sum_{k=1}^4Int_{9j}\left(\Omega_k\right)+Int_{9j}\left(\Omega^{out}\right)\right)
\end{align*}
where $\Omega_1,\Omega_2,\Omega_3,\Omega_4$ and $\Omega^{out}$ defined in the proof of Lemma \ref{lem:pPsi}. 

It follows from symmetry that $p_1=-2p_{10}$, $p_3=-p_4$, and $p_7=p_8$. Indeed, using the expression of $\I_{ij}$, we have $\KK\cdot K^r_1(x,y)=-2\KK\cdot K^r_{10}(x,y)$,  $\KK\cdot K^r_3(x,y)=-\KK\cdot K^r_4(y,x)$ and $\KK\cdot K_7^r(x,y)=\KK\cdot K_8^r(y,x)$. Since $\Theta^r$ is a radial function, then $f^r(x,y)=f^r(y,x)$. It follows from \eqref{J-I} that $p_1=-2p_{10}$, $p_3=-p_4$, and $p_7=p_8$.

It is easy to see that $Int_{21}\left(\R2\right)=0$ since $\I_{21}=0$. We have
\begin{align*}
\begin{split}
    Int_{91}\left(\R2 \right)=&\ \frac{1}{r^2}\int_{\R2}\frac{128  (x^4-y^4) |\psi_r|^2f^r(x,y)}{(1+|\psi_r|^2)^4}dxdy\\
    =&\ \frac{1}{r^2}\int_{\R2}  \left(\frac{128  (x^4-y^4)}{(1+|\psi_r|^2)^3}-\frac{128  (x^4-y^4)}{(1+|\psi_r|^2)^4}\right)f^r(x,y)dxdy.  
\end{split}
\end{align*}
Now let $p(x,y,r)=128  (x^4-y^4)  $,  it holds that $p(1,1,1)=p(-1,-1,1)=0$. Then we shall refine the proof of  Lemma \ref{lem:pPsi-2} to compute $Int_{91}\left(\R2\right)$. Since $f^r=1$ on $\Omega_1$, similar to \eqref{IOmega1-2nd} with $k=4,l=3,4$ and $p(x,y,r)=128  (x^4-y^4)  $, we have
\begin{align}\label{p1-Omega1}
\begin{split}
    Int_{91}(\Omega_1)=&\ \frac{\pi[\Delta_{x,y}p(\mathbf{1})-2 p_x(\mathbf{1})-2p_y(\mathbf{1})]}{128\cdot3\cdot2\cdot1}r^{-4}+O(r^{-8})  = O(r^{-8}).
\end{split} 
\end{align}
Note that $f^r=-1$ on $\Omega_2$,  we also have
\begin{align}\label{p1-Omega2}
\begin{split}
    Int_{91}(\Omega_2)=&\ -\frac{\pi[\Delta_{x,y}p(-\mathbf{1})+2 p_x(-\mathbf{1})+2p_y(-\mathbf{1})]}{128\cdot3\cdot2\cdot1}r^{-4}+O(r^{-8})
    =O(r^{-8}).
\end{split}
\end{align}
Note that $|f^r|\leq1$, similar to \eqref{IOmega3} and \eqref{IOmegaout}, we have
\begin{align}\label{p1-else}
    \left|Int_{91}\left(\Omega_3\right)\right|+ \left|Int_{91}\left(\Omega_4\right)\right|+\left|Int_{91}\left(\Omega^{out}\right)\right|= O(r^{-\frac{13}{2}}).
\end{align}
From \eqref{p1-Omega1}, \eqref{p1-Omega2} and \eqref{p1-else}, we get $p_1=-2p_{10}= O(r^{-\frac{13}{2}})$. Using the same method, we have
\begin{align}\label{pj-cp}
  p_3= -p_4=O(r^{-\frac92}),\quad p_5=O(r^{-\frac92}),\quad p_6=O(r^{-\frac92}).
\end{align}

We have
\begin{align*}
    Int_{22}(\R2)=&\ \int_{\R2}\frac{128(x^2+y^2)|\psi_r|^2f^r(x,y)}{(1+|\psi_r|^2)^4}dxdy\\
    =&\ \int_{\R2}  \left(\frac{128  (x^2+y^2)}{(1+|\psi_r|^2)^3}-\frac{128  (x^2+y^2)}{(1+|\psi_r|^2)^4}\right)f^r(x,y)dxdy.  
\end{align*}
We shall refine the proof of Lemma \ref{lem:pPsi-2} to compute $Int_{22}\left(\R2\right)$.  As before, similar to \eqref{IOmega1} with $k=2, l=3,4$ and $p(x,y,r)=128(x^2+y^2)$, we have
\begin{align}\label{p2-Omega1}
    \begin{split}
      Int_{22}\left(\Omega_1\right) =&\ \frac{\pi p(\mathbf{1})}{8\cdot3\cdot2\cdot1}+\frac{\pi[\Delta_{x,y}p(\mathbf{1})-2 p_x(\mathbf{1})-2p_y(\mathbf{1})+2p(\mathbf{1})]}{128\cdot3\cdot2\cdot1}r^{-4}+O(r^{-6})\\
      =&\ \frac{16\pi}{3}+O(r^{-6}).
    \end{split}
\end{align}
Note that $f^r=-1$ on $\Omega_2$, similar to \eqref{IOmega2}, we have
\begin{align}\label{p2-Omega2}
    \begin{split}
      Int_{22}\left(\Omega_2\right) =&\ -\frac{\pi p(-\mathbf{1})}{8\cdot3\cdot2\cdot1}-\frac{\pi[\Delta_{x,y}p(-\mathbf{1})+2 p_x(-\mathbf{1})+2p_y(-\mathbf{1})+2p(-\mathbf{1})]}{128\cdot3\cdot2\cdot1}r^{-4}+O(r^{-6})\\
      =&\ -\frac{16\pi}{3}+O(r^{-6}).
    \end{split}
\end{align}
Similar to \eqref{IOmega3} and \eqref{IOmegaout}, we have
\begin{align}\label{p2-else}
    \left|Int_{22}\left(\Omega_3\right)\right|+ \left|Int_{22}\left(\Omega_4\right)\right|+\left|Int_{22}\left(\Omega^{out}\right)\right|= O(r^{-6}).
\end{align}
From \eqref{p2-Omega1},\eqref{p2-Omega2} and
\eqref{p2-else}, we have $Int_{22}\left(\R2\right)=O(r^{-6})$. Similarly, we can also get $Int_{92}\left(\R2\right)=O(r^{-6})$. Then $p_2=O(r^{-6})$. Using the same method, we have 
\begin{align*}
     p_7=p_8=-\frac{16\pi}{3}r^{-4}+O(r^{-6}),\quad p_9=O(r^{-6}).
\end{align*}

Next, we want to solve the linear system $\J\vec{c}=\vec{p}$. To that end, we shall use the expression of $\J$ in \eqref{J-switch} after some row  and column switching. It is reduced to solving each block independently. For instance, we solve $A_1(c_1,c_{10},c_6)^T=(p_1,p_{10},p_6)^T$. Using Cramer's rule and \eqref{pj-cp},
\begin{align} 
    c_1=&\frac{1}{\det A_1}\begin{vmatrix}p_1& \J_{1,10}&\J_{16}\\ p_{10}&\J_{10,10}&\J_{10,6}\\p_6&\J_{6,10}&\J_{6,6}\end{vmatrix}\approx{r^4}\begin{vmatrix}O(r^{-\frac{13}{2}})&1&O(r^{-\frac92}) \\ O(r^{-\frac{13}{2}})&1&O(r^{-\frac{9}{2}})\\ O(r^{-\frac92})&O(r^{-\frac92})&1\end{vmatrix}=O(r^{-\frac52}).
\end{align}
Using the same method, we have $c_{10}=O(r^{-\frac52})$ and $c_6=O(r^{-\frac92})$. 
Moreover, one can verify that
\begin{align}
\begin{split}\label{c1-c10}
    c_1-2c_{10}=\frac{1}{\det A_1}\begin{vmatrix}p_1& \J_{1,10}+2\J_{11}&\J_{16}\\ p_{10}&\J_{10,10}+2\J_{1,10}&\J_{10,6}\\p_6&\J_{6,10}+2\J_{16}&\J_{6,6}\end{vmatrix}\approx&\ r^4\begin{vmatrix}O(r^{-\frac{13}{2}})&O(r^{-6})&O(r^{-\frac92}) \\ O(r^{-\frac{13}{2}})&r^{-4}&O(r^{-\frac{9}{2}})\\ O(r^{-\frac92})&O(r^{-\frac92})&1\end{vmatrix}\\
    =&\ O(r^{-\frac{13}{2}}).
    \end{split}
\end{align}
To solve $A_2(c_2,c_9,c_5)^T=(p_2,p_9,p_5)^T$, we use Cramer's rule to obtain
\begin{align} 
    c_2=&\frac{1}{\det A_2}\begin{vmatrix}p_2& \J_{29}&\J_{25}\\ p_{9}&\J_{99}&\J_{95}\\p_5&\J_{59}&\J_{55}\end{vmatrix}\approx r^4\begin{vmatrix}O(r^{-6})&1&O(r^{-\frac92})\\ O(r^{-6})&1&O(r^{-\frac92})\\O(-\frac92)&O(r^{-\frac92})&1\end{vmatrix}=O(r^{-2})
\end{align}
and $c_9=O(r^{-2})$, $c_5=O(r^{-\frac92})$. Furthermore,  similar to the approach of \eqref{c1-c10}, we can derive $c_2+2c_{9}=O(r^{-6})$.

To solve $A_3(c_3,c_8)^T=(p_3,p_8)^T$, 
\begin{align} 
    \begin{split}
     & c_3=\frac{1}{\det A_3}\begin{vmatrix}p_3&\J_{38}\\p_8&\J_{88}\end{vmatrix}\approx \begin{vmatrix}O(r^{-\frac92})&r^{-2}\\r^{-4}&1\end{vmatrix}=O(r^{-\frac92}),    \\
     &  c_8=\frac{1}{\det A_3}\begin{vmatrix}\J_{33}&p_3\\\J_{83}&p_8\end{vmatrix}=-\frac14r^{-4}+O(r^{-6}).
    \end{split}
\end{align}  
Using $p_3=-p_4$, $p_7=p_8$ and $\J_{38}=-\J_{47}$, we obtain $c_{4}=-c_3$ and $c_7=c_8$.
\end{proof}

\begin{proof}[\bf Proof of Lemma \ref{lem:nablafK}]
Making the following expansion
\begin{align}\label{nab_fK}
    \int_{\R2}|\nabla (f^r\KK)|^2=\int_{\R2}|\nabla f^r||\KK|^2+(f^r)^2|\nabla \KK|^2+2f^r\partial_\alpha f^r\KK_i\cdot \partial^\alpha \KK_i.
\end{align}
Integration by parts
\begin{align*}
    \int_{\R2}2f^r\partial_\alpha f^r\KK_i\cdot \partial^\alpha \KK_i=\int_{\R2}\partial_\alpha (f^r)^2\KK_i\cdot \partial^\alpha \KK_i=-\int_{\R2}(f^r)^2[\KK\cdot \Delta \KK+|\nabla \KK|^2].
\end{align*}
Since $\KK\in \L[\mathcal{S}(\Psi[\vec{\alpha}_r])]$, we have 
\begin{align}
    \int_{\R2}(f^r)^2\KK\cdot \Delta \KK=-\int_{\R2}|\nabla \mathcal{S}(\Psi[\vec{\alpha}_r])|^2|f^r\KK|^2.
\end{align}
Inserting the above two equations back to \eqref{nab_fK}, we get  
\begin{align} 
\begin{split}
    \int_{\R2}|\nabla (f^r\KK)|^2-|\nabla \mathcal{S}(\Psi[\vec{\alpha}_r])|^2|f^r\KK|^2
    =-\int_{\R2}|\nabla f^r|^2|\KK|^2.
    \end{split}
\end{align}
Next, we want to compute the right-hand side. Recall that $\KK=2K_2^r-K_9^r$. Then
\begin{align} 
\begin{split}
    |\KK|^2=4|K_2^r|^2-4K_2^r\cdot K_9^r+|K_9^r|^2
    =&\ \frac{4|\psi_r|^2}{(1+|\psi_r|^2)^2}\left[4-\frac{8xy}{r^2}+\frac{(x^2+y^2)^2}{r^4}\right]\\
    =&\ \frac{4|\psi_r|^4}{r^4(1+|\psi_r|^2)^2}.
    \end{split}
\end{align}
Introduce the notation $A_1=\{(x,y):r\leq (x-r)^2+(y-r)^2<r^2\}$ and $A_2=\{(x,y):r<(x+r)^2+(y+r)^2<r^2\}$. 
Recall that $f^r$ is defined in \eqref{def:fr}. Then
\begin{align} 
    |\nabla f^r|^2= \frac{4}{|\log r|^2}\frac{1}{(x-r)^2+(y-r)^2}\chi_{A_1}+\frac{4}{|\log r|^2}\frac{1}{(x+r)^2+(y+r)^2}\chi_{A_2}
\end{align}
where $\chi_{A_1}$ is the characteristic function of set $A_1$.
Then
\begin{align} 
\begin{split}
  \int_{A_1}|\nabla f^r|^2|\KK|^2\leq&\ \frac{16}{|\log r|^2r^4}\int_{A_1}\frac{1}{(x-r)^2+(y-r)^2}\frac{|\psi_r|^4}{(1+|\psi_r|^2)^2}\\
=&\ \frac{16}{|\log r|^2r^4}\int_{A_1}\frac{[(x+r)^2+(y+r)^2]|\psi_r|^2}{(1+|\psi_r|^2)^2}\\
\leq &\ \frac{16}{|\log r|^2r^4}C|\log r|=\frac{C}{|\log r|r^4}, 
\end{split}
\end{align}
for some uniform constant $C$. Similar estimate holds on $A_2$. 
Therefore
\begin{align} 
    \int_{\R2}|\nabla f^r|^2|\KK|^2=\int_{A_1}|\nabla f^r|^2|\KK|^2+\int_{A_2}|\nabla f^r|^2|\KK|^2=O(\frac{1}{|\log r|r^4}).
\end{align}
On the other hand, applying Corollary \ref{cor:int-psi}, one gets
\begin{align*}
    \int_{\R2}|\nabla\mathcal{S}(\Psi[\vec{\alpha}_r])|^2|f^r\KK|^2
    =&\int_{\Omega_1\cup \Omega_2}|\nabla\mathcal{S}(\Psi[\vec{\alpha}_r])|^2|\KK|^2+\int_{\Omega_3\cup\Omega_4}|\nabla\mathcal{S}(\Psi[\vec{\alpha}_r])|^2|f^r\KK|^2\\
    =&\ 4\J_{22}-4\J_{29}+\J_{99}+O(r^{-6})=\frac{64\pi}{3}r^{-4}+O(r^{-6}).
\end{align*} 
\end{proof}

\section*{Acknowledgement}
The research of L. Sun and J. Wei are partially supported by the National Key R\&D Program of China (No. 2022YFA1005601 and No. 2022YFA1005602). J. Wei is also partially supported by NSERC of Canada. The research of  B. Deng is partially supported by China Postdoctoral Science Foundation No. 2022M722474 and  NSFC  No. 11971137.

\bibliographystyle{plainnat}
\bibliography{ref.bib}
\end{document}